\documentclass[12pt,leqno]{amsart}
\usepackage[latin1]{inputenc}
\usepackage[T1]{fontenc}
\usepackage [pdftex] {graphicx}
\usepackage{color}
\usepackage[francais]{babel}
\usepackage{geometry}
\usepackage{amsmath,amsfonts,amssymb}

\newtheorem{thm}{Théorème}[section]
\newtheorem{cor}[thm]{Corollaire}
\newtheorem{lem}[thm]{Lemme}
\newtheorem{prop}[thm]{Proposition}
\newtheorem{defn}[thm]{Définition}
\newtheorem{rem}[thm]{Remarque}

\numberwithin{equation}{section}

\newcommand{\R}{{\mathbb{R}}}

\newcommand{\C}{{\mathbb{C}}}
\newcommand{\D}{{\mathrm{D}}}
\newcommand{\pd}{{\mathrm{d}}}
\newcommand{\h}{{\mathbb{H}}}

\newcommand{\RE}{{\mathrm{Re}}}
\newcommand{\dzi}{{\dfrac{\partial}{\partial{z_i}}}}
\newcommand{\dzj}{{\dfrac{\partial}{\partial{z_j}}}}

\newcommand{\dzn}{{\dfrac{\partial}{\partial{z_n}}}}
\newcommand{\dznb}{{\dfrac{\partial}{\partial{\overline{z}_n}}}}

\begin{document}
\title{Analyticité des applications CR dans des variétés presque complexes.}
\begin{abstract}
Nous étudions l'analyticité des applications CR entre deux hypersurfaces dans des variétés presque complexes. Nous démontrons l'analyticité d'une telle application lorsque la structure presque complexe d'arrivée est une déformation d'une structure modèle et lorsque l'hypersurface d'arrivée a une fonction définissante particulière. La preuve utilise la méthode de prolongation des systèmes d'équations aux dérivées partielles ainsi que la théorie des systèmes complets. 
\end{abstract}
\maketitle

\section{Introduction}
Le travail fondamental de M. Gromov (\cite{G85}) établit des résultats de référence en géométrie symplectique. Les courbes pseudo-holomorphes constituent un des outils essentiels de ce travail. L'analyse dans les variétés presque complexes se situe donc au c\oe ur des préoccupations en géométrie symplectique. L'étude des disques pseudo-holomorphes et de leur régularité au bord (\cite{IR04}), ainsi que des propriétés des domaines strictement pseudoconvexes (\cite{EG89}, \cite{M91}) ont déjà permis d'obtenir des résultats en géométrie symplectique et en géométrie de contact.\

Depuis le théorème de Newlander-Nirenberg (\cite{NN57}), nous savons qu'une structure presque complexe n'est pas nécessairement intégrable. De nombreux travaux ont naturellement cherché à savoir quels résultats demeurent dans le cadre presque complexe et quels sont les résultats qui font défaut. Ainsi, l'équivalence entre la positivité de la forme de Lévi et l'existence de  fonctions d'exhaustion strictement plurisousharmoniques permet, comme dans le cadre complexe, de caractériser les domaines pseudoconvexes dans les variétés presque complexes (\cite{DS08}). En revanche, le théorème de Wong-Rosay n'est plus vérifié dans le cadre presque complexe (\cite{L06}). Pour démontrer les résultats qui demeurent, il est nécessaire d'adapter les techniques de l'analyse complexe.\

Nous nous intéressons ici à l'analyticité des applications CR entre deux variétés presque complexes. Dans le cadre complexe, un tel résultat concernant l'analyticité d'une application CR entre deux variétés a été établi pour la première fois par S. Pinchuk dans \cite{P75} et H. Lewy dans \cite{L77} en utilisant un principe de réflexion : si les variétés $M$ et $M'$ sont strictement pseudoconvexes et analytiques réelles, un difféomorphisme CR de classe $\mathcal{C^1}$ entre $M$ et $M'$ s'étend en une application holomorphe. Le principe de réflexion est aussi utilisé par F. Forstneric dans \cite{F89} et par J. Faran dans \cite{F90}. \

De plus, M. S. Baouendi, J. Jacobowitz et F. Treves (\cite{BJT}) ont démontré l'analyticité des difféomorphismes $\mathcal{C}^\infty$ entre deux variétés pseudoconvexes CR dont l'une des deux est essentiellement finie. Une hypothèse supplémentaire de prolongation est nécessaire afin d'autoriser l'utilisation du ``edge of the wedge theorem''.\ 

Enfin, C. K. Han obtient des résultats d'analyticité en utilisant la théorie de prolongation et des systèmes complets, dans \cite{H97} et \cite{H83} notamment. Il s'agit d'obtenir un système complet dont la fonction CR $f$ est solution en différençiant les équations de Cauchy-Riemann tangentielles vérifiées par cette fonction $f$.\

Dans le cas presque complexe, nous pouvons utiliser la théorie de la prolongation et démontrer le Théorème \ref{prop2}, qui constitue le résultat principal de cet article. Il s'agit d'une généralisation au cas presque complexe d'un théorème dû à C. K. Han (\cite{H97}), pour certaines déformations des structures modèles vérifiant la condition $(*)$ (voir la section 4.1 pour la définition précise).
\begin{thm} \label{prop2}\
Soient $\Gamma$ et $\Gamma'$ deux hypersurfaces dans $\C^n$ telles que
$\Gamma=\{z\in\C^n,\,\rho(z)=0\}$, où $\rho$ est une fonction analytique réelle telle que telle que $\rho(z)=\RE(z_n)+|z'|²+o(|z|)$ et $\Gamma'=\{w\in\C^n,\,\rho'(w)=0\}$, où $\rho'$ est une fonction analytique réelle telle que $\rho'(w)=\RE(w_n)+|w'|²+o(|w|)$. Soient $J$ et $J'$ deux structures presque complexes analytiques réelles sur $\C^n$, avec $J'$ une structure vérifiant la condition $(*)$.\\ 
Soit $U$ une boule ouverte centrée en $0$ dans $\C^n$. Soit $g \,:\, \Gamma \cap U \rightarrow \Gamma'$ une application CR de classe $\mathcal{C}^4$. On suppose que $g$ un difféomorphisme local en 0 et que $g(0)=0$. Alors $g$ est uniquement déterminée par ses jets d'ordre 2 en un point, et $g$ est analytique réelle.
\end{thm}
On a donc :
\begin{cor}
Soit $\Gamma$ une hypersurface définie sur $\C^n$ par $\Gamma=\{w\in\C^n,\,\rho(w)=0\}$, où $\rho$ est une fonction analytique réelle telle que $\rho(w)=\RE(w_n)+|w'|²+o(|w|)$. Soit $J$ une structure presque complexe sur $\C^n$ vérifiant la condition $(*)$. Alors, le groupe des automorphismes CR de $(\Gamma,J)$ est de dimension finie et $dim_\R Aut(\Gamma,J)\leqslant (2n)^3$.
\end{cor}
Remarquons que cette estimée sur la dimension n'est pas précise.\\
Nous signalons aussi que la démonstration du Théorème \ref{prop2} s'adapte dans le cas où les hypersurfaces et les structures presque complexes sont lisses. On obtient alors que l'application est uniquement déterminée par ses jets d'ordre 2 en un point.\\

La preuve du Théorème \ref{prop2} constitue un premier pas vers la démonstration d'un résultat plus général concernant l'analyticité des applications CR entre deux hypersurfaces strictement pseudoconvexes dans des variétés presque complexes analytiques réelles.\\

La théorie des systèmes complets (\cite{H08}) permet de démontrer que certaines applications sont analytiques réelles. Il suffit d'exprimer toutes les dérivées d'un certain ordre $k$ d'une fonction $f$ vérifiant un système d'équations différentielles comme des fonctions analytiques réelles des dérivées d'ordre inférieur ou égal à $k$-$1$ de cette fonction pour conclure à l'analyticité de la fonction $f$. (\textit{cf.} 1.4) \

Les structures modèles, introduites dans \cite{CGS05} et étudiées dans \cite{CGS06}, jouent un rôle de référence dans l'étude des propriétés des domaines strictement pseudoconvexes.  Elles apparaissent naturellement comme limite d'un processus de dilatation anisotrope des coordonnées et constituent l'analogue presque complexe de la boule dans le cadre complexe. En effet, dans le cadre complexe, le théorème de
Wong-Rosay (\cite{W77}, \cite{R79}) stipule que, à biholomorphisme près, le seul domaine strictement pseudoconvexe à groupe d'automorphisme non compact est la boule unité. Il s'agit un résultat local puisqu'un domaine strictement pseudoconvexe qui admet une orbite d'accumulation en un point de stricte pseudoconvexité de sa frontière est biholomorphe à la boule unité (\cite{GKK02}). Ce résultat fait défaut dans le cadre presque complexe. En effet, les structures modèles admettent des orbites d'accumulation en un point de stricte $J$-pseudoconvexité de leur frontière. De plus, un domaine qui admet une orbite d'accumulation en un point de stricte $J$-pseudoconvexité de sa frontière est biholomorphe à une variété modèle (\cite{CGS06}, \cite{L06}).\

Les structures modèles ont notamment été utilisées par H. Gaussier et A. Sukhov (\cite{CGS06}) pour démontrer qu'un difféomorphisme lisse entre deux domaines $D$ et $D'$ lisses et relativement compacts dans des variétés réelles s'étend à la frontière à condition que $D$ admette une structure presque complexe lisse sur $\overline{D}$ telle que $(D,J)$ soit strictement pseudoconvexe et que la structure $f_\star(J)$ admette une extension lisse sur $\overline{D'}$ avec $(D',f_\star(J))$ strictement pseudoconvexe. \

On peut noter que toute variété presque complexe strictement pseudoconvexe est une petite déformation d'une variété modèle. Ce résultat classique est obtenu en utilisant une technique de dilatation anisotrope des coordonnées. Cette méthode constitue un analogue presque complexe de la technique des dilatations introduite par S. Pinchuk (\cite{P91}) dans le cadre complexe.\\

Les structures presque complexes vérifiant la condition $(*)$ sont particulières puisque le défaut d'intégrabilité de $T^{1,0}_J \C^n$ est porté par une seule direction. Les structures vérifiant la condition $(*)$ constituent une généralisation des structures modèles. Dans le cas de la dimension réelle 4, la seule structure modèle est la structure standard. En revanche, la matrice complexe représentant une structure presque complexe vérifiant la condition $(*)$ est de la forme suivante : \begin{eqnarray*}J_\C=\begin{pmatrix}
i&0&0&0\\
0&-i&0&0\\
a&b&i+c&d\\
\overline{b}&\overline{a}&\overline{d}&-i+\overline{c}\\                                                                                                                                                                                                                                                                                                                                                                                                                                                                                                                                                                                                                                                                                                                                                                                                                                                                                                                                                                                                                                                                                                                                                                                                                                                                                                                                                                                                                                      \end{pmatrix}.\end{eqnarray*}
Il existe donc des structures non standard qui vérifient la condition $(*)$. Il est toujours possible, en dimension réelle 4, d'obtenir une forme ``normale`` pour une structure presque complexe (\cite{S94}). La matrice complexe d'une telle structure est une matrice diagonale par blocs :
\begin{eqnarray*}J_\C=\begin{pmatrix}
J_1(z)&0\\
0&J_2(z)\\                                                                                                                                                                                                                                                                                                                                                                                                                                                                                                                                                                                                                                                                                                                                                                                                                                                                                                                                                                                                                                                                                                                                                                                                                                                                                                                                                                                                                                      \end{pmatrix}.\end{eqnarray*}
Nous verrons que, même en dimension 2, il n'est pas toujours possible d'obtenir, après un difféomorphisme local, une structure presque complexe vérifiant la condition $(*)$.\\


Dans le cas d'une application CR entre $(\Gamma,J)$ et $(\Gamma',J')$ vérifiant les hypothèses du Théorème \ref{prop2}, les équations de Cauchy-Riemann tangentielles  sont très proches des équations dans le cas complexe. La prolongation du système CR afin d'obtenir un système complet est essentiellement une adaptation de celle réalisée par C. K. Han dans \cite{H97}.\\

L'article est organisé comme suit : dans la première partie nous définissons les notions de base concernant les variétés presque complexes. La deuxième partie est consacrée à la démonstration du Théorème \ref{prop2} dans le cas particulier où l'hypersurface $\Gamma$ est définie par $\Gamma=\{z\in\C^n,\RE(z_n)+|z'|^2=0\}$, (où $z'=(z_1,\cdots,z_{n-1})$), et lorsque les structures presque complexes sont des structures modèles. Nous établissons dans le premier paragraphe le système d'équations de Cauchy-Riemann tangentielles vérifiées par la fonction $f$. Le deuxième paragraphe, très calculatoire, est consacré à la démonstration d'une proposition clé : la Proposition \ref{l4}. Il s'agit d'exprimer certaines dérivées de la fonction $f$ comme des fonctions analytiques réelles d'autres dérivées de la fonction $\overline{f}$. A l'aide de ce résultat, nous parvenons, dans le troisième paragraphe, à expliciter un système complet vérifié par la fonction $f$ et ainsi terminer la preuve. La troisième partie contient la démonstration du théorème \ref{prop2}. Nous définissons dans le premier paragraphe les structures presque complexes vérifiant la condition $(*)$. Dans le deuxième paragraphe, nous montrons que le système d'équations de Cauchy-Riemann tangentielles vérifiées par la fonction $g$ permet d'appliquer la même méthode que dans la démonstration du Théorème \ref{prop1} pour obtenir l'analyticité de la fonction $g$. Nous donnons dans le troisième paragraphe une interprétation géométrique pour les structures presque complexes vérifiant la condition $(*)$.
\tableofcontents

\section{Préliminaires}
\subsection{Variétés presque complexes, Applications CR }
Soit $M$ une variété différentielle de dimension réelle $2n$. Une\textbf{ structure presque complexe} $J$ sur $M$ est la donnée d'un isomorphisme de fibrés vectoriels différentiable, $J \,:\, T(M) \rightarrow T(M)  $ vérifiant  $J^2=-I$. On appelle \textbf{variété presque complexe} une variété différentielle $M$ munie d'une structure presque complexe $J$.\\
Une structure presque complexe $J$ est dite \textbf{intégrable} lorsque $(M,J)$ est une variété complexe.\\
Soient $(M,J)$ et $(M',J')$ deux variétés presque complexes. Une application $f\,:\,M \rightarrow M'$ de classe $\mathcal{C}^1$ est dite \textbf{$(J,J')$-holomorphe}, ou \textbf{pseudo-holomorphe}, si 
\begin{align}\label{ph}
\forall z \in M, \, \pd f_z\circ J_z=J'_{f(z)} \circ \pd f_z. 
\end{align}
Soit $\Gamma$ une sous-variété de $(M,J)$. Nous noterons  $T_{\C}\Gamma$ le complexifié de l'espace tangent de $\Gamma$: $ T_{\C}\Gamma=\C \otimes T\Gamma $. On définit l'espace tangent $J$-holomorphe de $\Gamma$ par $H^{1,0}\Gamma=\{Z \in T_{\C}\Gamma, JZ=iZ \}$. Il intervient dans la définition d'une application CR entre deux variétés presque complexes : 
\begin{defn}
Soient $(M,J)$ et $(N,J')$ deux variétés presque complexes. Soient $\Gamma$ et $\Gamma'$ deux sous-variétés de $M$ et $N$. Une application $\mathcal{C}^1$, $f : \Gamma\rightarrow \Gamma'$ est dite CR si $f_\star\{H^{1,0}\Gamma\} \subset H^{1,0}\Gamma'$.
\end{defn}
\begin{rem}\label{rq1}
Soit $f\,:\,(M,J)\mapsto(M',J')$ une bijection pseudo-holomorphe. Soit $\Gamma$ une sous-variété de $M$. Alors l'application $\tilde{f}=f_{|\Gamma}$ définie sur $\Gamma$ et à valeurs dans $f(\Gamma)=\Gamma'$ est CR.\\
Localement, soient $\Omega$ et $\Omega'$ des ouverts de $M$ et $M'$, tels que $\Gamma \cap \Omega $ est non vide, soit $f\,:\,(\Omega,J)\mapsto(\Omega',J')$ une bijection pseudo-holomorphe. Alors l'application $\tilde{f}=f_{|\Gamma\cap\Omega}$ définie sur $\Gamma\cap\Omega$ et à valeurs dans $f(\Gamma\cap\Omega)$ est CR.
\end{rem}

Si $\varphi$ est une 1-forme sur $M$, alors $J^\star\varphi$ est la forme définie sur l'espace tangent de $M$ par, $(J^\star\varphi) X=\varphi(JX)$, pour $X\in TM$. Le crochet de Lie de deux champs de vecteurs $X$ et $Y$ est le champ de vecteurs $[X,Y]$ tel que pour toute fonction $f$ de classe $\mathcal{C}^\infty$ sur $M$, on ait : $[X,Y]f=X(Yf)-Y(Xf)$.
Soit $\Gamma$ une hypersurface réelle lisse de $M$ définie par $\Gamma=\{r=0\}$. Soit $p\in\Gamma$. Le fibré tangent $J$-holomorphe de $\Gamma$ est défini par $H_p^J\Gamma=T_p\Gamma\cap JT_p\Gamma$. C'est sur cet espace qu'est définie la forme de Lévi : 
\begin{defn}\label{pseudoconvexe}
\begin{enumerate}
 \item La \textbf{forme de Lévi} de $\Gamma$ en $p$ est l'application définie sur $H_p^J\Gamma$ par $\mathcal{L}^J_\Gamma(X_p)=J^\star \pd r [X,JX]_p$, où $X$ est un champ de vecteurs de $H^J\Gamma$ tel que $X(p)=X_p$. (La définition ne dépend pas du choix d'un tel $X$).
\item Une variété presque complexe $(M,J)$ est dite \textbf{strictement $J$-pseudoconvexe} en $p$ si, pour tout $X_p$ tel que $X_p\in H_p^J\Gamma$, $\mathcal{L}^J_\Gamma(X_p)>0$.
L'hypersurface $\Gamma$ est dite \textbf{strictement $J$-pseudoconvexe} si elle l'est en tout point.
\end{enumerate}
\end{defn}

\subsection{Structures modèles}
Nous définissons maintenant les structures modèles. Nous utilisons les notations de H. Gaussier et A. Sukhov dans \cite{CGS06}. 
\begin{defn} \label{modèle}
Une structure presque complexe $J$ sur $\C^n$ est dite \textbf{structure modèle} si $J(z)=J_{st}+L(z)$, où $L$ est une matrice $L=(L_{j,k})_{1 \leqslant j,k \leqslant 2n} $ telle que 
\begin{align*}
 L_{j,k}&=0 \text{ si } 1\leqslant j \leqslant 2n-2, \,1 \leqslant k \leqslant 2n,\\
 L_{j,k}&=\sum_{l=1}^{n-1}(a_l^{j,k}z_l+\overline{a}_l^{j,k}\overline{z}_l) ,\, a_l^{j,k}\in \C, \text{ si } j=2n-1,2n \text{ et } k=1,\ldots , 2n-2.
\end{align*}
La complexification $J_\C$ d'une structure modèle s'écrit comme une matrice complexe $2n\times 2n$  :
\begin{eqnarray*}J_{\C}=\begin{pmatrix}
i&0&0&0&\ldots&0&0\\
0&-i&0&0&\ldots&0&0\\
0&0&i&0&\ldots&0&0\\
0&0&0&-i&\ldots&0&0\\
\ldots&\ldots&\ldots&\ldots&\ldots&\ldots&\ldots\\
0&\tilde{L}_{2n-1,2}(z,\overline{z})&0&\tilde{L}_{2n-1,4}(z,\overline{z})&\ldots&i&0\\
\tilde{L}_{2n,1}(z,\overline{z})&0&\tilde{L}_{2n,3}(z,\overline{z})&0&\ldots&0&-i\\
\end{pmatrix},\end{eqnarray*}
avec $\tilde{L}_{2n,2i-1}(z,\overline{z})=\sum_{l=1}^{n-1}(\alpha_l^i z_l+\beta_l^i\overline{z}_l)$, où  $\alpha_l^i,\, \beta_l^i \in \C$ et $\tilde{L}_{2n,2i-1}=\overline{\tilde{L}_{2n-1,2i}}$.\\
Soit $J$ une structure modèle sur $\C^n$, et $D=\{z\in \C^n,\,Re(z_n)+P(z',\overline{z'})=0\}$, avec $P$ un polynôme homogène du second degré sur $\C^{n-1}$ à valeurs réelles. Le couple $(D,J)$ est dit \textbf{domaine modèle} si $D$ est strictement $J$-pseudoconvexe au voisinage de l'origine.
\end{defn}
Remarquons qu'une structure modèle est nécessairement analytique réelle.\\
Nous étudions maintenant les structures modèles et l'espace tangent $J$-holomorphe pour l'hypersurface $\Gamma$ définie par $\Gamma=\partial\h=\{z\in\C^n,\RE(z_n)+|z'|^2=0\}$, où $z'=(z_1,\cdots,z_{n-1})$ et $\h$ est le demi-plan de Siegel défini par $\h=\{z\in\C^n,\RE(z_n)+|z'|^2<0\}$.\\
Soit $J_{mod}$ une structure modèle sur $\C^n$. Soit $H^{1,0}\partial\h=T\partial\h \cap JT\partial\h$ l'espace tangent $J$-holomorphe . Les $(n-1)$ champs de vecteurs \begin{equation}\label{champs}
L_j=\dzi+\alpha_j(z)\dzn+\beta_j(z)\dznb,\; j=1,\ldots , n-1,
\end{equation}
 forment une base de $H^{1,0}\partial\h$, avec
\begin{equation}\label{betai}\beta_j(z) =-\dfrac{i}{2}\tilde{L}_{2n,2j-1}(z)=-\dfrac{i}{2}\sum_{l=1}^{n-1}(\alpha_l^j z_l+\beta_l^j\overline{z}_l):=\sum_{l=1}^{n-1}(a_l^j z_l+b_l^j\overline{z}_l)
 \end{equation}
et,
\begin{equation}\label{alphai}
\alpha_j(z)=2(\dfrac{i}{4}\tilde{L}_{2n,2j-1}(z)-\overline{z}_j)=-\sum_{l=1}^{n-1}(a_l^j z_l+b_l^j\overline{z}_l)-2\overline{z}_j
\end{equation}
 Nous définissons un champ de vecteurs $T$ comme étant la projection du champ de vecteurs $[L_1,\overline{L_1}]$ dans $T_{\C}\partial\h / H^{1,0}\partial\h \bigoplus\overline{H^{1,0}\partial\h}$. Nous avons alors \[T_{\C}\partial\h=H^{1,0}\partial\h\bigoplus\overline{H^{1,0}\partial\h}\bigoplus<T>.\] Le calcul du crochet de Lie $[L_1,\overline{L}_1]$ permet de calculer explicitement le champ $T$. Nous avons $T=i(\dzn-\dznb)$. Ainsi, $\{T,\, L_j,\,\overline{L}_j,\, j=1,\ldots , n-1\}$ est une base de $T_{\C} \partial\h$, le complexifié de l'espace tangent de $\partial\h$.\\
En particulier, nous remarquons que 
\begin{equation}\label{rho1}
[L_j,\overline{L}_k]=\gamma_{j,\overline{k}}T,\;\text{ avec }\gamma_{j,\overline{k}}=\left\{\begin{split}
-2i +\dfrac{1}{2}\left(\beta_j^j+\overline{\beta}_j^j\right)\text { si } j=k \\
\dfrac{1}{2}\left(\beta_k^j+\overline{\beta}_j^k\right) \text{ si } j\neq k.\end{split}\right.                           
\end{equation}
Ainsi,
\begin{equation}\label{rho}
\overline{L}_kL_j=L_j\overline{L}_k-[L_j,\overline{L}_k]=L_j\overline{L}_k-\gamma_{j,\overline{k}}T.
\end{equation}
On remarque en particulier que pour tout $j=1,\ldots,n-1$, $\gamma_{j,\overline{j}}\neq0$.
Nous avons aussi $[\overline{L}_k,L_j]=\gamma_{\overline{k},j}T$, avec $\gamma_{\overline{k},j}=-\gamma_{j,\overline{k}}$.\\
Il faut noter que $\gamma_{j,\overline{k}}$ est constant. Enfin, nous remarquons que les champs $T$ et $L_k$, $(k=1,\ldots,n-1)$ commutent, ainsi que les champs $T$ et $\overline{L}_k$, $(k=1,\ldots,n-1)$.

\subsection{Systèmes complets}
La théorie des systèmes complets permet de conclure à l'analyticité des applications vérifiant un tel système. Dans ce paragraphe, nous définissons les systèmes complets. Nous utilisons les notations et nous rappellons les résultats de \cite{H08}.\

Soient $U$ et $V$ des ouverts de $\R^m$ et $\R^n$. Soit $f=(f_1,\ldots,f_n):\,U\rightarrow \,V$ une application de classe $\mathcal{C}^k$, vérifiant un système d'équations différentielles d'ordre $q$ ($q\leq k$), pour $x=(x_1,\ldots,x_m)\in U$ :
\begin{eqnarray}\label{complet'}
 \Delta_p(x,\D^\beta f, |\beta|\leq q)=0,\,p=1,\ldots,l,
\end{eqnarray}
où les applications $ \Delta_p$ sont lisses.
\begin{defn}
On dit que $f$ vérifie un système complet d'ordre $k$ lorsque toutes les dérivées partielles de $f_j$, $j=1,\ldots, n$, d'ordre $k$ peuvent être exprimées comme des fonctions lisses des derivées de $f_1,\ldots,f_n$ d'ordre inférieur à $k$ : \\
Pour tout $j=1,\ldots,n$, pour tout multi-indice $\alpha$ tel que $|\alpha|=k$, il existe $H_j^\alpha$ lisse telle que :  
\begin{equation}\label{complet}
\D ^\alpha f_j=H_j^\alpha(\D^\beta f,\,|\beta|<k)
\end{equation}
\end{defn} 
D'après \cite{H08}, nous avons : 
\begin{prop}\label{analytique}
Soit $f$ une application de classe $\mathcal{C}^k$ vérifiant un système complet (\ref{complet}) d'ordre $k$. Alors,
\begin{enumerate}
\item $f$ est uniquement déterminée par son jet d'ordre $(k-1)$ en un point, et $f$ est de classe $\mathcal{C}^\infty$.
\item Si de plus, les applications $ \Delta_p$ et $H_j^\alpha$ sont analytiques réelles (pour $p=1,\ldots,l$, $j=1,\ldots,n$ et $|\alpha|=k$ ), alors l'application $f$ est aussi analytique réelle.
\end{enumerate}
\end{prop}
La proposition \ref{analytique} a été utilisée par C. K. Han, en particulier dans \cite{H83} et \cite{H97}. Dans le cas où $k_0=1$, le théorème suivant constitue le cas particulier du Théorème \ref{prop1} pour des structures complexes.
\begin{thm} \label{thmH97}(C.K. Han \cite{H97})\\
Soit $M$ une variété analytique réelle CR, Levi non dégénérée, de dimension $2m+1$. Soit $\{L_1,\ldots,L_m\}$ une base du fibré $\mathcal{V}$ définissant la structure CR. Soit $N$ une sous-variété analytique réelle de $\C^{n+1}$, ($n\geqslant m$) définie par $r(z,\overline{z})=0$ (où $r$ est normalisée).\\
Soit $f \,:\, M\rightarrow N$ une application CR telle que, pour un certain entier $k_0$, les vecteurs $\{L^\alpha f,\,|\alpha|\leqslant k_0\}$ et $(0,\ldots,0,1)$ engendrent $\C^{n+1}$.\\ 
Si $f$ est de classe $\mathcal{C}^{2k_0+1}$, alors $f$ est analytique réelle.
\end{thm}
Pour démontrer que les fonctions considérées vérifient un système complet, C. K. Han utilise la théorie de la prolongation. Prolonger un système d'équations aux dérivées partielles consiste à différencier ce système un certain nombre de fois. Génériquement, en différençiant le système d'équations vérifiées par la fonction autant de fois que nécessaire, il est possible d'inverser le système prolongé et d'exprimer les dérivées partielles d'un certain ordre $k$ comme des fonctions lisses, ou analytiques réelles, des dérivées partielles d'ordre inférieur ou égal à $k-1$.


\section{Démonstration dans le cas des structures modèles.}
Nous démontrons dans cette partie un cas particulier du Théorème \ref{prop2} lorsque l'hypersurface $\Gamma$ est définie par $\Gamma=\partial\h=\{z\in\C^n,\RE(z_n)+|z'|^2=0\}$, (où $z'=(z_1,\cdots,z_{n-1})$), et lorsque les structures presque complexes sont des structures modèles.
\begin{thm} \label{prop1}\
Soient  $J_{mod}$ et $J'_{mod}$ deux structures modèles sur $\C^n$.\\
Soit $U$ une boule ouverte centrée en $0$ dans $\C^n$. Soit $f \,:\, \partial\h\cap U \rightarrow \partial\h$ une application CR de classe $\mathcal{C}^4$. On suppose que $f$ est un difféomorphisme local en 0 et que $f(0)=0$. Alors $f$ est uniquement déterminée par ses jets d'ordre 2 en un point, et $f$ est analytique réelle.
\end{thm}
Pour prouver le Théorème \ref{prop1}, nous utiliserons la proposition \ref{analytique} à propos des systèmes complets pour obtenir l'analyticité de $f$. Il suffira de montrer que toutes ses dérivées d'ordre 3 s'expriment de façon analytique réelle en fonction de ses dérivées d'ordre inférieur ou égal à 2.

\subsection{Système d'équations vérifiées par $f$.}
Nous cherchons maintenant à écrire une condition nécessaire et suffisante pour que $f$ soit CR de $(\partial\h\cap U,J_{mod})$ dans $(\partial\h, J'_{mod})$ : \\
Soit $(L_1(z),\ldots,L_{n-1}(z))$ une base de $H^{1,0}_{J_{mod}}\partial\h$ et $(Z_1(w),\ldots,Z_{n-1}(w))$ une base de $H^{1,0}_{J'_{mod}}\partial\h$.\\
Si $f$ est CR, pour chaque $L_p(z)$, nous avons : $f_\star\{L_p(z)\} \in H^{1,0}_{J'_{mod}}\partial\h$. \\
Calculons d'abord  $f_\star\{L_p(z)\}$ : \\
\begin{equation}\label{CR1}
f_\star\{L_p(z)\}=\sum_{j=1}^n \left( L_p(z)f_j(z) \dfrac{\partial}{\partial{w_j}} + L_p(z)\overline{f}_j(z)\dfrac{\partial}{\partial{\overline{w}_j}}\right).
\end{equation}
De plus, $H^{1,0}_{J'_{mod}}\partial\h=<Z_1(w),\ldots, Z_{n-1}(w)>$, avec $w=f(z)$.\\
$f_\star\{L_p(z)\}$ s'écrit donc : 
\begin{align}
f_\star\{L_p(z)\}&=a_1(w) Z_1(w)+\ldots + a_{n-1}(w) Z_{n-1}(w).\label{CR1bis}
\end{align}
D'après (\ref{champs}), nous avons $Z_j(w)=\dfrac{\partial}{\partial{w_j}}+\beta_j(w)\dfrac{\partial}{\partial{\overline{w}_n}}+\alpha_j(w)\dfrac{\partial}{\partial{w_n}}$. Ainsi, dans les égalités (\ref{CR1}) et (\ref{CR1bis}), les termes en $\dfrac{\partial}{\partial{w_j}}$  et $\dfrac{\partial}{\partial{\overline{w_j}}}$ doivent être égaux. On obtient donc 
\begin{align}
  a_j(w)=L_p(z) f_j(z), \text{ pour }j=1,\ldots,n-1,\notag\\
L_p(z) \overline{f}_j(z)=0, \text{ pour }j, p=1,\ldots ,n-1.\label{CR2}
\end{align}
De plus, 
\begin{align*}
f_\star\{L_p(z)\}=&\sum_{j=1}^{n-1}  L_p(z)f_j Z_j(w)  \\
 & +\left(L_p(z)f_n(z)-\sum_{j=1}^{n-1}\alpha_{j}(f(z))L_p(z) f_j(z)\right)\dfrac{\partial}{\partial{w_n}}\\
&+ \left(L_p(z)\overline{f}_n(z)-\sum_{j=1}^{n-1}\beta_j(f(z))L_p(z) f_j(z)\right)\dfrac{\partial}{\partial{\overline{w}_n}}.
\end{align*}
Nous avons donc, pour $p=1,\ldots,n-1$ :
\begin{align}
L_p(z)f_n(z)-\sum_{j=1}^{n-1}\alpha_{j}(f(z))L_p(z) f_j(z)=0,\label{CR3}\\
L_p(z)\overline{f}_n(z)-\sum_{j=1}^{n-1}\beta_j(f(z))L_p(z) f_j(z)=0.\label{CR4}
\end{align}
Avec les égalités (\ref{CR2}), (\ref{CR3}) et (\ref{CR4}), nous obtenons :\\
$f=(f_1,\ldots,f_n)$ est une application CR de $\partial\h \cap U$ dans $\partial\h$ si et seulement si les égalités suivantes sont vérifiées sur $\partial\h \cap U$ : 
\begin{equation}\label{3.1}
L_p \overline{f}_j=0, \text{ pour }j, p=1,\ldots ,n-1 ,
\end{equation}
\begin{equation}\label{3.1bis}
L_p \overline{f}_n=\sum_{j=1}^{n-1}\beta_j(f)L_p f_j,  \text{ pour }p=1,\ldots ,n-1 ,
\end{equation}
\begin{equation}\label{3.1ter}
 L_p f_n=\sum_{j=1}^{n-1}\alpha_{j}(f)L_p f_j, \text{ pour }p=1,\ldots ,n-1 ,
\end{equation}
\begin{equation}\label{3.2}
\dfrac{f_n+\overline{f}_n}{2}+\sum_{j=1}^{n-1}f_j\overline{f}_j=0.
\end{equation}
où (\ref{3.2}) est l'écriture de $ f(\partial\h\cap U)\subset \partial\h$.\\
Nous remarquons que si $J'_{mod}=J_{st}$, alors $\beta_j(z)=0$ pour $j=1,\ldots,n-1$. L'équation (\ref{3.1bis}) devient donc $L_p \overline{f}_n=0$, pour $p=1,\ldots,n-1$, et nous retrouvons les équations dans le cas complexe de \cite{H97}.\\


\subsection{Prolongation du système CR}
Soit $\mathcal{C}_{p,q}$ l'espace des fonctions analytiques réelles en les variables
$$\{\,T^tL^\alpha f_j : t+|\alpha|\leqslant p\, , t\leqslant q ,\, j=1,\ldots,n-1,\; L^\alpha f_n : |\alpha|\leqslant p\}.$$ 
L'introduction de cet espace, qui apparaît naturellement dans les calculs, permet d'alléger les notations. Nous noterons $\overline{\mathcal{C}}_{p,q}$ l'espace des fonctions analytiques réelles en les variables 
$$\{\,T^t\overline{L}^\alpha \overline{f}_j : t+|\alpha|\leqslant p\, , t\leqslant q ,\, j=1,\ldots,n-1,\; \overline{L}^\alpha \overline{f}_n : |\alpha|\leqslant p\}.$$ 
Nous allons prouver dans cette partie la proposition suivante : 
\begin{prop}\label{l4} Pour $p=1,\ldots ,n, $
\begin{eqnarray*}
(i)&\text{Pour tous }t,\,\alpha\text{ tels que }t+|\alpha|\leqslant 3,\, T^t L^\alpha f_p \in \overline{\mathcal{C}}_{1+t,1}.\\
(ii)&\text{Pour tous }t,\,\alpha\text{ tels que }t+|\alpha|\leqslant  4, \text{ et }t\leqslant 1, \,T^t L^\alpha f_p \in \overline{\mathcal{C}}_{2,1}.
\end{eqnarray*}
\end{prop}

C'est l'étape clé dans la démonstration du Théorème \ref{prop1} puisque l'utilisation de l'assertion $(i)$ de ce résultat et du conjugué de l'assertion $(ii)$ permettront d'obtenir, avec des calculs supplémentaires, un système complet vérifié par $f$.
La preuve de cette proposition est adaptée de la démonstration de \cite{H97}. Elle est scindée en quatre propositions, les Propositions \ref{l1}, \ref{l2}, \ref{l3} et \ref{l5}, correspondant chacune aux valeurs de $t$ égales à $0$, $1$, $2$ et $3$.\\

Nous commençons par démontrer le lemme suivant : 
\begin{lem}\label{l0}
La matrice $(L_k f_j)_{k,j=1,\ldots,n-1}$ est inversible en 0.
\end{lem}
\noindent \textsc{Preuve du Lemme \ref{l0}} : Puisque $f$ est un difféomorphisme CR local en 0, nous avons $<f_*(L_1)(0),\ldots,  f_*(L_{n-1})(0)>=H^{1,0}_J\partial\h$.Or, 
\begin{align*}
f_\star(L_k)(0)=\sum_{j=1}^n \left( L_p(0)f_j(0) \dfrac{\partial}{\partial{w_j}} + L_p(0)\overline{f}_j(0)\dfrac{\partial}{\partial{\overline{w}_j}}\right).
\end{align*}
De plus, d'après (\ref{3.1}), (\ref{3.1bis}) et (\ref{3.1ter}), on a, 
\begin{align*}
L_k(0)& \overline{f}_j=0, \text{ pour }j=1,\ldots ,n-1 ,\\
L_k(0)& \overline{f}_n(0)=\sum_{j=1}^{n-1}\beta_j(f(0))L_k(0) f_j(0)=0,  \text{ car } f(0)=0 \text{ et d'après (\ref{betai}), }\beta_j(0)=0,\\
L_k(0)&f_n(0)=\sum_{j=1}^{n-1}\alpha_{j}(f(0))L_p(0) f_j(0)=0, \text{ car } f(0)=0 \text{et d'après (\ref{alphai}), }\alpha_j(0)=0.
\end{align*}
Ainsi, les $(n-1)$ vecteurs $f_*(L_k)(0)=\sum_{j=1}^{n-1} L_p(0)f_j(0) \dfrac{\partial}{\partial{w_j}}$ engendrent $H^{1,0}_J\partial\h$ qui est de dimension $(n-1)$. La matrice $(L_k f_j)_{k,j=1,\ldots,n-1}$ est donc inversible en 0.$\square$\\
\\

Appliquons le champ $\overline{L}_k$ à l'égalité (\ref{3.2}). Avec (\ref{3.1}), il vient : 
\begin{equation}\label{3.3}
\dfrac{\overline{L}_k f_n+\overline{L}_k \overline{f}_n}{2}+\sum_{j=1}^{n-1}f_j\overline{L}_k \overline{f}_j=0
\end{equation}

Considérons alors le système linéaire d'équations (\ref{3.3}) (pour $k=1,\ldots,  n-1$) et (\ref{3.2}) d'inconnues $(f_1,\ldots,f_n)$. La matrice  $ \begin{pmatrix}
\overline{L}_1 \overline{f}_1&\ldots&\overline{L}_1 \overline{f}_{n-1}&0\\
\vdots&\vdots&\vdots&\vdots\\
\overline{L}_{n-1} \overline{f}_1&\ldots&\overline{L}_{n-1} \overline{f}_{n-1}&0\\
\overline{f}_1&\ldots&\overline{f}_{n-1}&\dfrac{1}{2}
\end{pmatrix}$ étant inversible en $0$, elle est inversible sur un voisinage de $0$. Nous pouvons résoudre ce système pour $(f_1,\ldots,f_n)$ sur un voisinage de 0 en fonction de $\{\overline{f}_j, \overline{L}_k \overline{f}_j; j=1,\ldots ,n, \overline{L}_k f_n ; k=1,\ldots ,n-1\}$. Ainsi,
\begin{equation*}
 f_p=H_p(\overline{f}_j, \overline{L}_k \overline{f}_j; j=1,\ldots ,n, \overline{L}_k f_n ; k=1,\ldots ,n-1),\;\text{pour } p=1,\ldots ,n,
\end{equation*}
où $H_p$ est une fonction analytique des termes à l'intérieur de la parenthèse.\\
En prenant le conjugué de l'égalité (\ref{3.1bis}), nous observons que $\overline{L}_k f_n$  peut s'écrire en fonction de $\{\overline{f}_j,\;\overline{L}_k \overline{f}_j,\; j=1,\ldots ,n-1\}$. Ainsi, 
\begin{equation}\label{3.4bis-}
 f_p=H_p(\overline{f}_j, \overline{L}_k \overline{f}_j; j=1,\ldots ,n,\;k=1,\ldots ,n-1),
\end{equation}
où $H_p$ est une autre fonction analytique des termes à l'intérieur de la parenthèse. Nous pouvons réecrire l'égalité (\ref{3.4bis-}) sous la forme :
\begin{equation}\label{3.4bis}
 f_p=H_p(\overline{f},\overline{L}\overline{f}).
\end{equation}Cette écriture est le point de départ de la démonstration de la Proposition \ref{l4}.\

Démontrons maintenant la Proposition \ref{l1}, qui constitue la première étape de la démonstration de la Proposition \ref{l4} :
\begin{prop}\label{l1} Pour  $p=1,\ldots, n$, pour $k=1,\ldots, n-1$,\\
\begin{eqnarray}
(i)&   L_k f_p \in \overline{\mathcal{C}}_{1,1} \label{l1.1}\\ 
(ii)&L_k(\overline{\mathcal{C}}_{1,1})\subset \overline{\mathcal{C}}_{1,1} \nonumber\\
&\text{En particulier, }  L^{\alpha} f_p \in \overline{\mathcal{C}}_{1,1} \text{ pour tout multi-indice } \alpha \text{ avec } |\alpha|\leqslant 4\label{l1.3}
\end{eqnarray}

\end{prop}

\noindent\textsc{Preuve de la Proposition \ref{l1}} : \textit{(i)} Fixons $k_0$ et appliquons  $L_{k_0}$ à l'égalité (\ref{3.4bis}) : 
\begin{equation*}
L_{k_0}f_p=L_{k_0}H_p(\overline{f},\overline{L}\overline{f}).
\end{equation*}
La fonction $H_p$ étant analytique réelle (en les $n²$ variables $\overline{f}_j, \overline{L}_k \overline{f}_j$, $j=1,\ldots,n$, $k=1,\ldots,n-1$) , nous pouvons l'écrire sous la forme \[H_p(\overline{f},\overline{L}\overline{f})=\sum_{\alpha} a_\alpha^p \overline{f}_1^{\alpha_1}\ldots \overline{f}_n^{\alpha_n}\overline{L}_1\overline{f}_1^{\alpha_{n+1}}\ldots\overline{L}_{n-1}\overline{f}_1^{\alpha_{n²-n+1}}\ldots \overline{L}_{n-1}\overline{f}_n^{\alpha_{n²}},\]
où la série $\sum_\alpha a_\alpha^p r_1^{\alpha_1}\ldots r_{n²}^{\alpha_{n²}}$ converge pour $r_j<r_0$, $1 \leqslant j \leqslant n²$. Ainsi,
\[L_{k_0}f_p=\sum_{\alpha} a_\alpha^p L_{k_0}(\overline{f}_1^{\alpha_1}\ldots \overline{L}_{n-1}\overline{f}_n^{\alpha_{n²}}).\]
Lorsque le champ $L_{k_0}$ s'applique à l'un des facteurs du terme $\overline{f}_1^{\alpha_1}\ldots \overline{L}_{n-1}\overline{f}_n^{\alpha_{n²}}$, il faut distinguer quatre cas : 
\begin{itemize}
 \item Le champ $L_{k_0}$ s'applique à $\overline{f}_j,\,j\neq n$ : d'après (\ref{3.1}), nous avons $L_{k_0}\overline{f}_j=0$.
 \item Lorsque le champ $L_{k_0}$ s'applique à $\overline{f}_n$, écrivons l'égalité (\ref{3.1bis}) et remplaçons $\beta_j(f)$ par sa valeur donnée dans (\ref{betai})  :
\begin{eqnarray}
 L_{k_0}\overline{f}_n&=&\sum_{j=1}^{n-1}\beta_j(f)L_{k_0} f_j \nonumber\\
&=&\sum_{j=1}^{n-1}\sum_{l=1}^{n-1}\left(a_j^lf_l+b_j^l \overline{f}_l\right)L_{k_0} f_j \label{Lkbfn}.
\end{eqnarray}
Ce cas se produit pour les multi-indices $\alpha$ tels que $\alpha_n\neq 0$.
 \item Le champ $L_{k_0}$ s'applique à $\overline{L}_m\overline{f}_j,\,j\neq n$. Nous utilisons d'abord (\ref{rho}): 
\begin{align}
L_{k_0}\overline{L}_m\overline{f}_j&=\overline{L}_mL_{k_0}\overline{f}_j-[\overline{L}_m,L_{k_0}]\overline{f}_j. \notag
\end{align}
D'après (\ref{3.1}), nous avons $L_{k_0}\overline{f}_j=0$. Nous pouvons de plus utiliser l'égalité (\ref{rho1}) et remplacer $[\overline{L}_m,L_{k_0}]$ par sa valeur.
\begin{align}
L_{k_0}\overline{L}_m\overline{f}_j&=0-\gamma_{\overline{m},k_0}T\overline{f}_j  \label{Lk0Lmfj}
\end{align}
Ce cas se produit pour les multi-indices $\alpha$ tels que $\alpha_{nm+j}\neq 0$, $m,j=1,\ldots,n-1$.
 \item Lorsque le champ $L_{k_0}$ s'applique à $\overline{L}_m\overline{f}_n$, nous utilisons le conjugué de l'égalité (\ref{3.1ter}) : \\
\begin{eqnarray}
L_{k_0}\overline{L_m}\overline{f_n}&=& L_{k_0}\left(\sum_{j=1}^{n-1}\overline{\alpha}_{j}(f)\overline{L}_m \overline{f}_j\right) .\nonumber
\end{eqnarray}
Nous remplaçons $\overline{\alpha}_{j}(f)$ par sa valeur, donnée par le conjugué de l'égalité (\ref{alphai}) :
\begin{eqnarray*}
L_{k_0}\overline{L_m}\overline{f_n}&= &L_{k_0}\left(\sum_{j=1}^{n-1}-\left(\left(\sum_{l=1}^{n-1}\overline{a}_j^l \overline{f}_l+\overline{b}_j^l f_l\right) -2f_j\right)\overline{L}_m \overline{f}_j\right)\\
&=&\sum_{j=1}^{n-1}\left(\left(\sum_{l=1}^{n-1}-\overline{b}_j^l L_{k_0}f_l -2L_{k_0}f_j\right)\overline{L}_m \overline{f}_j\right.\nonumber\\
&&\hspace{0,7cm}\left.-\left(\left(\sum_{l=1}^{n-1}\overline{a}_j^l \overline{f}_l+\overline{b}_j^l f_l \right)-2f_j\right)L_{k_0}\overline{L}_m \overline{f}_j\right).
\end{eqnarray*}
Nous avons vu en (\ref{Lk0Lmfj}) que $L_{k_0}\overline{L}_m\overline{f}_j=-\gamma_{\overline{m},k_0}T\overline{f}_j$. Nous avons donc : 
\begin{eqnarray}
L_{k_0}\overline{L_m}\overline{f_n}&=&\sum_{j=1}^{n-1}\left(\left(\sum_{l=1}^{n-1}-\overline{b}_j^l L_{k_0}f_l -2L_{k_0}f_j\right)\overline{L}_m \overline{f}_j \right.\nonumber\\
&&\hspace{0,7cm}\left.+ \left(\left(\sum_{l=1}^{n-1}\overline{a}_j^l \overline{f}_l+\overline{b}_j^l f_l\right) -2f_j\right)\gamma_{\overline{m},k_0} T\overline{f}_j\right).\label{LmLkbfnb}
\end{eqnarray}
Ce cas se produit pour les multi-indices $\alpha$ tels que $\alpha_{nm}\neq 0$, $m=1,\ldots,n$.
\end{itemize}
Nous pouvons donc écrire l'égalité suivante : 
\begin{align*}
 L_{k_0}f_p&=h_{k_0}^p+\sum_{k=1}^{n-1} \varphi_{k_0}^{p,k}\left(\overline{f}_l,f_l,\overline{L}_q\overline{f}_l,\,q=1,\ldots,n-1,\,l=1,\ldots,n\right)L_{k_0}f_k, \text{ pour } p=1,\ldots,n,\\
\end{align*}
avec
\begin{eqnarray*}
h_{k_0}^p &=&\sum_{j,m=1,\ldots,n-1}\sum_{\alpha\,|\, \alpha_{nm+j}\neq0} \psi^{1,p,k_0}_{\alpha,m,j}+
 \sum_{m=1,\ldots,n}\sum_{\alpha\,|\, \alpha_{mn}\neq0} \psi^{2,p,k_0}_{\alpha,m}
\end{eqnarray*}
et, 
\begin{eqnarray*}
\varphi_{k_0}^{p,k}&=&\sum_{\alpha\,|\, \alpha_n\neq0} \psi^{3,p,k_0}_{\alpha,n,k} +\sum_{m=1}^{n_-1}\sum_{\alpha\,|\, \alpha_{mn}\neq0}\psi^4_{\alpha,m,n,k,k_0} +\sum_{\alpha\,|\, \alpha_{(k+1)n}\neq0}\psi^{5,p,k_0}_{\alpha,k} ,
\end{eqnarray*}
où,
\begin{eqnarray*}
\psi^{1,p,k_0}_{\alpha,m,j}&=&-\gamma_{\overline{m},k_0}T\overline{f}_j\alpha_{nm+j}a_\alpha^p \overline{f}_1^{\alpha_1}\ldots\overline{f}_n^{\alpha_n-1}\ldots\overline{L}_m\overline{f}_j^{\alpha_{nm+j}-1}\ldots\overline{L}_{n-1}\overline{f}_n^{\alpha_{n²}},\\
\psi^{2,p,k_0}_{\alpha,m}&=&\alpha_{mn}a_\alpha^p((\sum_{l=1}^{n-1}\overline{a}_j^l \overline{f}_l+\overline{b}_j^l f_l )-2f_j)\gamma_{\overline{m},k_0} T\overline{f}_j \overline{f}_1^{\alpha_1}\ldots \overline{L}_m\overline{f}_{n}^{\alpha_{mn}-1}\ldots \overline{L}_{n-1}\overline{f}_n^{\alpha_{n²}},\\
\psi^{3,p,k_0}_{\alpha,n,k}&=&\alpha_na_\alpha^p \left(\sum_{l=1}^{n-1} a_k^l f_l+b_k^l \overline{f}_l\right)\overline{f}_1^{\alpha_1}\ldots \overline{f}_n^{\alpha_n-1}\ldots \overline{L}_{n-1}\overline{f}_n^{\alpha_{n²}},\\
\psi^{4,p,k_0}_{\alpha,m,n,k}&=& (-2\overline{L}_m \overline{f}_k+\sum_{j=1}^{n-1}-\overline{b}_j^k \overline{L}_m \overline{f}_j)\alpha_{mn}a_\alpha^p \overline{f}_1^{\alpha_1}\ldots \overline{L}_m\overline{f}_{n}^{\alpha_{mn}-1}\ldots \overline{L}_{n-1}\overline{f}_n^{\alpha_{n²}},\\
\psi^{5,p,k_0}_{\alpha,k}&=&-2\overline{L}_k \overline{f}_k\alpha_{(k+1)n}a_\alpha^p \overline{f}_1^{\alpha_1}\ldots \overline{L}_k\overline{f}_{n}^{\alpha_{(k+1)n}-1}\ldots \overline{L}_{n-1}\overline{f}_n^{\alpha_{n²}}.\\
\end{eqnarray*}

Ainsi, en notant $A_{k_0}=\begin{pmatrix}
1-\varphi_{k_0}^{1,1}&\ldots&-\varphi_{k_0}^{1,n}\\
\vdots&&\vdots\\
-\varphi_{k_0}^{n,1}&\ldots&1-\varphi_{k_0}^{n,n}\\
\end{pmatrix}$, on a $\begin{pmatrix}
A_{k_0}
\end{pmatrix}\begin{pmatrix}
L_{k_0}f_1\\
\vdots\\
L_{k_0}f_n\\
\end{pmatrix}=\begin{pmatrix}
h_{k_0}^1\\
\vdots\\
h_{k_0}^n\\
\end{pmatrix}.$\\
Nous remarquons que $h_{k_0}^p$ appartient à l'espace $\overline{\mathcal{C}}_{1,1}$, ainsi que $\varphi_{k_0}^{p,k}$.\\
Soit $\delta_0$ tel que si une matrice $B=(b_{i,j})_{1\leqslant i,j\leqslant n}$ vérifie $|b_{i,j}|<\delta_0$ pour $1\leqslant i,j\leqslant n$, alors la matrice  $A=I_n-B$ est inversible.\\
Nous savons que $\sum_\alpha \sum_{q=1,\ldots, n²} \alpha_q|a_\alpha^p| r_0^{\alpha_1}\ldots r_0^{\alpha_q-1}\ldots r_{0}^{\alpha_{n²}}=M<\infty$. \\
Soit $b=\sup\{2+\sum_{j=1}^{n-1}|\overline{b}_j^k |,\,k=1,\ldots,n-1\}$.\\
Soit $\epsilon=\min(r_0,\dfrac{\delta_0}{bM})$.
\begin{lem}\label{dilatation} \
\begin{enumerate}
\item Après dilatation, nous pouvons supposer que :  $\forall\, j,k=1,\ldots, n-1\,:\, |L_kf_j(0)|<\epsilon$.
\item La matrice $A_{k_0}(0)$ est alors inversible.
\end{enumerate}
\end{lem}

\noindent\textsc{Preuve du Lemme \ref{dilatation} : } (i) : Soit $z=(z',z_n)\in \C^{n-1}\times\C $.\\
Soit $\lambda_\delta : \partial\h  \rightarrow  \partial\h$ l'application définie par $\lambda_\delta(z',z_n)  =(\sqrt{\delta}z',\delta z_n)$.\\
Soit $f^\delta=\lambda_\delta \circ f$.
L'application $\lambda_\delta $ est $(J,J')$-holomorphe sur $\C^n$ et $\lambda_\delta (\partial\h)=\partial\h$, donc d'après la remarque \ref{rq1}), l'application $\lambda_\delta $ étant un automorphisme CR de $\partial\h$. L'application $f^\delta$ vérifie donc aussi les hypothèses du Théorème \ref{prop1}. Ainsi, $f^\delta$ est analytique réelle si et seulement si $f$ l'est.\\
De plus, si $j\neq n$ : $L_kf_j^\delta=L_k(\lambda_\delta \circ f)_j=\sqrt{\delta}L_kf_j$.\\
De même, $L_kf_n^\delta(0)=\delta L_kf_n$.
Ainsi, pour $\delta$ assez petit, nous avons bien $|L_kf^\delta_j(0)|<\epsilon$.\\

(ii) : Nous avons la majoration suivante :  
\begin{eqnarray}|\varphi_{k_0}^{p,k}(0)| &\leqslant&\sum_{m=1}^{n_-1}\sum_{\alpha\,|\, \alpha_{mn}\neq0}|\psi^4_{\alpha,m,n,k,k_0}| +\sum_{\alpha\,|\, \alpha_{(k+1)n}\neq0}|\psi^{5,p,k_0}_{\alpha,k}|\label{M1}\end{eqnarray}
Or,
\begin{eqnarray*}
|\psi^4_{\alpha,m,n,k,k_0}| &\leqslant &2|\overline{L}_k \overline{f}_k(0)|\alpha_{(k+1)n}|a_\alpha^p| |\overline{f}_1(0)|^{\alpha_1}\ldots |\overline{L}_k\overline{f}_{n}(0)|^{\alpha_{(k+1)n}-1}\ldots |\overline{L}_{n-1}\overline{f}_n(0)|^{\alpha_{n²}}
\end{eqnarray*}
Ainsi, puisque $|L_kf_j(0)|<\epsilon=\min(r_0,\dfrac{\delta_0}{bM})$, nous avons : 
\begin{eqnarray}
|\psi^4_{\alpha,m,n,k,k_0}| &\leqslant &2\epsilon\alpha_{(k+1)n}|a_\alpha^p| r_0^{\alpha_1}\ldots r_0^{\alpha_{(k+1)n}-1}\ldots r_0^{\alpha_{n²}}\nonumber\\
&\leqslant &\epsilon b\alpha_{(k+1)n}|a_\alpha^p| r_0^{\alpha_1}\ldots r_0^{\alpha_{(k+1)n}-1}\ldots r_0^{\alpha_{n²}}.\label{M2}
\end{eqnarray}
De même, nous obtenons
\begin{eqnarray}
|\psi^{5,p,k_0}_{\alpha,k}|&\leqslant&(2\epsilon+\sum_{j=1}^{n-1}|\overline{b}_j^k|\epsilon)\alpha_{(m+1)n}|a_\alpha^p| r_0^{\alpha_1}\ldots r_0^{\alpha_{(m+1)n}-1}\ldots r_0^{\alpha_{n²}}\nonumber\\
&\leqslant&\epsilon b\alpha_{(m+1)n}|a_\alpha^p| r_0^{\alpha_1}\ldots r_0^{\alpha_{(m+1)n}-1}\ldots r_0^{\alpha_{n²}}.\label{M3}
\end{eqnarray}
En replaçant les majorations obtenues en (\ref{M2}) et (\ref{M3}) dans l'inégalité (\ref{M1}), nous obtenons : 
\begin{eqnarray*}|\varphi_{k_0}^{p,k}(0)| &\leqslant  & \epsilon b\sum_{\alpha} \sum_{q=1}^{n^2}\alpha_{q}|a_\alpha^p| r_0^{\alpha_1}\ldots r_0^{\alpha_{q}-1}\ldots r_0^{\alpha_{n²}}\\
&\leqslant&  \epsilon bM \leqslant \delta_0.
\end{eqnarray*}
La matrice $A_{k_0}(0)$ est donc inversible.$\square$\\

La matrice $A_{k_0}$ est donc inversible sur un voisinage de 0. Ainsi,  $\begin{pmatrix}
L_{k_0}f_1\\
\vdots\\
L_{k_0}f_n\\
\end{pmatrix}=\begin{pmatrix}
A_{k_0}
\end{pmatrix}^{-1}\begin{pmatrix}
h_{k_0}^1\\
\vdots\\
h_{k_0}^n\\
\end{pmatrix}.$\\
On a $A_{k_0}=I_n-B$, avec $B=(\varphi_{k_0}^{p,k})_{p,k=1,\ldots,n}$. Rappellons que les $\varphi_{k_0}^{p,k}$ appartiennent à l'espace $\overline{\mathcal{C}}_{1,1}$. D'après la formule $(I_n-B)^{-1}=\sum_{p=0}^{+\infty}(-B)^p$, les coefficients de la matrice $(A_{k_0})^{-1}$ sont dans l'espace $\overline{\mathcal{C}}_{1,1}$. Puisque $h^p_{k_0} \in \overline{\mathcal{C}}_{1,1}$ pour $p=1,\ldots,n$, nous avons $L_{k_0}f_p \in \overline{\mathcal{C}}_{1,1}$, pour $p=1,\ldots, n$.

\vspace{0.5cm}\
\underline{\textit{ (ii)  Montrons que pour $m=1,\ldots ,n-1$; $L_m(\overline{\mathcal{C}}_{1,1})\subset \overline{\mathcal{C}}_{1,1}$.}}\\
Il suffit de voir que les termes $L_m \overline{f}_j$, $L_m \overline{L}_k\overline{f}_j$, $L_m T \overline{f}_j$, $L_m \overline{f}_n$, $L_m \overline{L}_k\overline{f}_n$ sont dans l'espace $\overline{\mathcal{C}}_{1,1}$ (pour $m,j,k=1,\ldots, n-1$) : 
\begin{itemize}
\item D'après (\ref{3.1}), nous avons $L_m \overline{f}_j =0$.
\item D'après l'égalité (\ref{Lk0Lmfj}), nous avons : \begin{eqnarray*}L_m \overline{L}_k \overline{f}_j  =-\gamma_{\overline{k},m} T \overline{f}_j  \in \overline{\mathcal{C}}_{1,1}.\end{eqnarray*}
\item Les champs $T$ et $L_m$ commutent, donc $L_m T \overline{f}_j = T L_m \overline{f}_j=0$.
\item D'après l'égalité (\ref{3.1bis}), nous avons : $L_m \overline{f}_n=\sum_{j=1}^{n-1}\beta_j(f)L_m f_j $.\\
Or, nous avons démontré au (i) que $L_m f_j \in \overline{\mathcal{C}}_{1,1} $, et d'après (\ref{betai}), nous avons $\beta_j(f)=\sum_{l=1}^{n-1}\left(a_j^l f_l+b_j^l \overline{f}_l \right)\in \overline{\mathcal{C}}_{1,1}$. Nous avons donc $L_m \overline{f}_n \in \overline{\mathcal{C}}_{1,1}$.
\item Pour le terme $L_m \overline{L}_k\overline{f}_n$, le calcul a déjà été fait en (\ref{LmLkbfnb}) :
\begin{eqnarray*} L_m \overline{L}_k\overline{f}_n=
\sum_{j=1}^{n-1}\left(\left(\sum_{l=1}^{n-1}-\overline{b}_j^l L_{m}f_l -2L_{m}f_j\right)\overline{L}_k \overline{f}_j + \gamma_{\overline{k},m} \left(\sum_{l=1,}^{n-1}\left(\overline{a}_j^l \overline{f}_l+\overline{b}_j^l f_l\right) -2f_j\right)T\overline{f}_j\right).
\end{eqnarray*}
Nous avons déjà démontré que $L_m f_l \in \overline{\mathcal{C}}_{1,1} $ et que $f_j\in \overline{\mathcal{C}}_{1,1}$. Ainsi, $L_m \overline{L}_k\overline{f}_n \in \overline{\mathcal{C}}_{1,1}$.
\end{itemize}
Par récurrence, nous obtenons immédiatemment que pour tout multi-indice $\alpha,\,|\alpha|\leqslant 4,\, L^\alpha f_p \in \overline{\mathcal{C}}_{1,1}$. Ceci termine la preuve de la Proposition \ref{l1}.$\square$\\

Le Lemme \ref{l00} sera utilisé régulièrement dans les prochains calculs : 
\begin{lem}\label{l00} Pour $m=1,\ldots ,n-1 $, pour $t=1,2,3$ :
\begin{eqnarray*}
\overline{L}_m ( \overline{\mathcal{C}}_{t,1}) \subset \overline{\mathcal{C}}_{1+t,1}.
\end{eqnarray*}
\end{lem}
\noindent\textsc{Preuve du Lemme \ref{l00}} : La démonstration se fait par récurrence sur $t$ : \\
\underline{Pour $t=1$} : Il suffit de voir que les termes $\overline{L}_m \overline{f}_j$, $\overline{L}_m \overline{L}_k\overline{f}_j$, $\overline{L}_m T \overline{f}_j$, $\overline{L}_m \overline{f}_n$, $\overline{L}_k \overline{L}_m\overline{f}_n$ sont dans $\overline{\mathcal{C}}_{2,1}$ (pour $m,j,k=1,\ldots, n-1$) : 
\begin{itemize}
\item Par définition, nous avons : $\overline{L}_m \overline{f}_j \in \overline{\mathcal{C}}_{1,1}$.
\item Par définition encore, nous avons : $\overline{L}_m \overline{L}_k \overline{f}_j \in \overline{\mathcal{C}}_{2,1}$.
\item  Puisque les champs $T$ et $\overline{L}_m$ commutent, nous avons : $\overline{L}_m T \overline{f}_j = T \overline{L}_m \overline{f}_j\in \overline{\mathcal{C}}_{2,1}$ par définition.
\item D'après le conjugué de l'égalité (\ref{3.1ter}), nous avons :
\begin{eqnarray*}
\overline{L}_m \overline{f}_n&=\sum_{j=1}^{n-1}\overline{\alpha}_j(f)\overline{L}_m \overline{f}_j. 
\end{eqnarray*}
Nous remplaçons $\overline{\alpha}_j(f)$ par sa valeur, donnée par le conjugué de (\ref{alphai}) : 
\begin{eqnarray}
\overline{L}_m \overline{f}_n&= -\sum_{j=1}^{n-1} \left(\sum_{l=1}^{n-1}\left(\overline{a}_j^l \overline{f}_l+\overline{b}_j^l f_l\right) +2f_j\right)\overline{L}_m \overline{f}_j \label{Lmbfnb}\\
&\in \overline{\mathcal{C}}_{1,1} \text{ car } f_j\in \overline{\mathcal{C}}_{1,1}.\notag
\end{eqnarray}
\item Appliquons maintenant le champ $\overline{L}_k$ à l'égalité (\ref{Lmbfnb}) :
\begin{align*} \overline{L}_k \overline{L}_m\overline{f}_n&=\overline{L}_k\left(\sum_{j=1}^{n-1}\overline{\alpha}_j(f)\overline{L}_m \overline{f}_j\right)\\
& =-\overline{L}_k\sum_{j=1}^{n-1}\left(\sum_{l=1}^{n-1}\left(\overline{a}_j^l \overline{f}_l+\overline{b}_j^l f_l\right)+2f_j\right)\overline{L}_m \overline{f}_j\\
& =-\sum_{j=1}^{n-1}\left(\left(\sum_{l=1}^{n-1}\left(\overline{a}_j^l \overline{f}_l+\overline{b}_j^l f_l\right)+2f_j\right)\overline{L}_k\overline{L}_m \overline{f}_j+\sum_{l=1}^{n-1}\overline{a}_j^l \overline{L}_k\overline{f}_l\overline{L}_m \overline{f}_j\right)\\
&\in \overline{\mathcal{C}}_{2,1} \text{ car }f_j\in \overline{\mathcal{C}}_{1,1}.
\end{align*}
\end{itemize}

\underline{$t-1\Rightarrow t$} : \\ 
Il suffit de voir que $\overline{L}_m T \overline{L}_{i_1}\ldots \overline{L}_{i_{t-1}}\overline{f}_j \in \overline{\mathcal{C}}_{1+t,1}$ (les autres termes sont donnés par l'hypothèse de récurrence ou par la définition).\\
Puisque les champs $T$ et $\overline{L}_m$ commutent, nous avons : 
\begin{eqnarray*}
\overline{L}_m T \overline{L}_{i_1}\ldots \overline{L}_{i_{t-1}}\overline{f}_j&=T\overline{L}_m  \overline{L}_{i_1}\ldots \overline{L}_{i_{t-1}}\overline{f}_j & \in \overline{\mathcal{C}}_{t+1,1}.
\end{eqnarray*}
$\square$\

Démontrons maintenant la Proposition \ref{l2}, qui constitue la deuxième étape de la démonstration de la Proposition \ref{l4} :
\begin{prop}\label{l2} Pour $p=1,\ldots ,n $, pour $k,m=1,\ldots ,n-1 $, : 
\begin{eqnarray*}
(i)&T f_p \in \overline{\mathcal{C}}_{2,1} \\
(ii)&T L_k f_p \in \overline{\mathcal{C}}_{2,1}\\
(iii)&T L_kL_m f_p \in \overline{\mathcal{C}}_{2,1}\\
(iv) &T L_kL_m L_j f_p \in \overline{\mathcal{C}}_{2,1}
\end{eqnarray*}
\end{prop}

\noindent\textsc{Preuve de la Proposition \ref{l2}} : \textit{(i) } Appliquons le champ $\overline{L}_1$ à l'égalité (\ref{l1.1}). D'après le Lemme \ref{l00}, nous avons $\overline{L}_1(\overline{\mathcal{C}}_{1,1})\subset \overline{\mathcal{C}}_{2,1}$. Ainsi, $\overline{L}_1 L_1 f_p\in \overline{\mathcal{C}}_{2,1}$.\\
Utilisons l'égalité (\ref{rho}) : 
\begin{eqnarray}
 \overline{L}_1 L_1 f_p&=& L_1\overline{L}_1 f_p-[L_1,\overline{L}_1] f_p\nonumber\\
&=& L_1\overline{L}_1 f_p -\gamma_{\overline{1},1}Tf_p\label{L1bL1fp}
\end{eqnarray}

\begin{itemize}
 \item  \underline{Si $p \neq n $ :} D'après (\ref{3.1}), nous avons $\overline{L}_1 f_p=0$. Ainsi, 
\begin{eqnarray}
\overline{L}_1 L_1 f_p&=  -\gamma_{\overline{1},1}Tf_p \label{L1bL1fj}
\end{eqnarray}
Ainsi, puisque $\gamma_{1,\overline{1}}\neq 0$,  $Tf_p =\dfrac{-1}{\gamma_{1,\overline{1}}}\overline{L}_1 L_1 f_p \in \overline{\mathcal{C}}_{2,1}$ (puisque $\gamma_{1,\overline{1}}$ est constant).
\item \underline{Si $p = n $ :} 
\begin{eqnarray*}
\overline{L}_1 L_1 f_n&=  L_1 \overline{L}_1 f_n -\gamma_{1,\overline{1}}Tf_n
\end{eqnarray*}
Or, d'après le conjugué de l'égalité (\ref{3.1bis}), nous avons :
\begin{eqnarray*}
L_1\overline{L}_1 f_n&=L_1\left(\sum_{j=1}^{n-1}\overline{\beta}_j(f)\overline{L}_1 \overline{f}_j\right)
\end{eqnarray*}
Nous remplaçons $\overline{\beta}_j(f)$ par sa valeur, donnée par le conjugué de (\ref{betai}) : 
\begin{eqnarray*}
L_1\overline{L}_1 f_n&=&L_1\left(\sum_{j=1}^{n-1}\sum_{l=1}^{n-1}\left(\overline{a}_j^l \overline{f}_l+\overline{b}_j^l f_l\right)\overline{L}_1 \overline{f}_j\right)\\
&=&\sum_{j=1}^{n-1}\sum_{l=1}^{n-1} \left(\overline{b}_j^l \left(L_1f_l\overline{L}_1 \overline{f}_j+ f_lL_1\overline{L}_1 \overline{f}_j+\overline{a}_j^l\overline{f}_lL_1\overline{L}_1 \overline{f}_j\right)\right)
\end{eqnarray*}
D'après le conjugué de l'égalité (\ref{L1bL1fj}), nous avons $L_1\overline{L}_1 \overline{f}_j=-\overline{\gamma_{1,\overline{1}}}T \overline{f}_j$. Ainsi, 
\begin{eqnarray}
L_1\overline{L}_1 f_n&=&\sum_{j=1}^{n-1}\sum_{l=1}^{n-1} \left(\overline{b}_j^l \left(L_1f_l\overline{L}_1 \overline{f}_j-\gamma_{1,\overline{1}}f_lT \overline{f}_j\right)-\gamma_{1,\overline{1}}\overline{a}_j^l\overline{f}_lT \overline{f}_j\right)\label{L11fn}\\
&\in& \overline{\mathcal{C}}_{2,1} \text{ car } L_1f_l \text{ et } f_l\in \overline{\mathcal{C}}_{2,1}\nonumber
\end{eqnarray}
Nous avons donc $\overline{L}_1 L_1 f_n-  L_1 \overline{L}_1 f_n \in \overline{\mathcal{C}}_{2,1}$.\\
Ainsi, $Tf_n =\dfrac{1}{\gamma_{1,\overline{1}}}\left(L_1 \overline{L}_1 f_n-\overline{L}_1 L_1 f_n\right)\in \overline{\mathcal{C}}_{2,1}$.\\
\end{itemize}\

\underline{\textit{(ii) Montrons que $T L_k f_p \in \overline{\mathcal{C}}_{2,1}$}}\\
Nous avons démontré en (\ref{l1.3}) que $L_kL_kf_p \in \overline{\mathcal{C}}_{1,1}$. Grâce au Lemme \ref{l00}, nous avons $\overline{L}_kL_kL_kf_p \in \overline{\mathcal{C}}_{2,1}$.\\
Nous chercons à calculer le terme $TL_kf_p$, pour cela, dans le terme $\overline{L}_kL_kL_kf_p$, nous faisons commuter deux fois les champs $\overline{L}_k$ et $L_k$ en utilisant l'égalité (\ref{rho}) : 
\begin{eqnarray*}
\overline{L}_kL_k L_k f_p&=& L_k \overline{L}_k L_k f_p - [\overline{L}_k, L_k]L_kf_p\\
&=&L_k L_k\overline{L}_k  f_p  - L_k [L_k,\overline{L}_k]  f_p - \gamma_{\overline{k},k}TL_kf_p\\
&=& L_k L_k\overline{L}_k  f_p - L_k(\gamma_{\overline{k},k}T ) f_p - \gamma_{\overline{k},k}TL_kf_p
\end{eqnarray*}
Le terme $\gamma_{\overline{k},k}$ étant constant, nous avons : 
\begin{eqnarray*}
\overline{L}_kL_k L_k f_p&=& L_k L_k\overline{L}_k  f_p - \gamma_{\overline{k},k}L_kT f_p- \gamma_{\overline{k},k}TL_kf_p
\end{eqnarray*}
Les champs $T$ et $L_k$ commutent. Nous avons donc : 
\begin{eqnarray*}
\overline{L}_kL_k L_k f_p&=&L_k L_k\overline{L}_k  f_p - \gamma_{\overline{k},k}TL_k f_p  - \gamma_{\overline{k},k}TL_kf_p\\
&=&L_k L_k\overline{L}_k  f_p - 2\gamma_{\overline{k},k}TL_k f_p 
\end{eqnarray*}
\begin{itemize}
 \item \underline{Si $p\neq n$ :} D'après (\ref{3.1}), nous avons $\overline{L}_k f_p=0$. Ainsi, 

\begin{eqnarray}
\overline{L}_kL_k L_k f_p&= & - 2\gamma_{\overline{k},k}TL_k f_p \label{TL1}
\end{eqnarray}
Nous obtenons $TL_k f_p = \dfrac{-1}{2\gamma_{k,\overline{k}}}\overline{L}_k L_k L_k f_p\in \overline{\mathcal{C}}_{2,1}$, puisque $\gamma_{\overline{k},k}$ est constant.\\

\item \underline{Si $p=n$ :} Reprenons le calcul en (\ref{TL1}):
\begin{align*}
\overline{L}_kL_k L_k f_n=L_k L_k\overline{L}_k  f_n - 2\gamma_{\overline{k},k}TL_k f_n.
\end{align*}
Le terme $L_k L_k\overline{L}_k  f_n$ est non nul. Pour le calculer, appliquons $L_k$ dans l'égalité (\ref{L11fn}) (en remplaçant 1 par $k$) : 
\begin{eqnarray*}
L_kL_k\overline{L}_k f_n
&=&L_k\left(\sum_{j=1}^{n-1}\sum_{l=1,}^{n-1} \left(\overline{b}_j^l( L_kf_l\overline{L}_k \overline{f}_j-\gamma_{k,\overline{k}}f_lT \overline{f}_j)\right)-\gamma_{k,\overline{k}}\overline{a}_j^l\overline{f}_lT \overline{f}_j\right)\\
&=&\sum_{j=1}^{n-1}\sum_{l=1}^{n-1}\left( \overline{b}_j^l (L_kL_kf_l\overline{L}_k \overline{f}_j+L_kf_lL_k\overline{L}_k \overline{f}_j\right.\\
&&\hspace{1,2cm}\left.-\gamma_{k,\overline{k}}( L_kf_lT \overline{f}_j+f_lL_kT \overline{f}_j))-\gamma_{k,\overline{k}}\overline{a}_j^lL_k(\overline{f}_lT \overline{f}_j)\right)
\end{eqnarray*}
Les champs $T$ et $L_k$ commutent puis le terme $L_k \overline{f}_j$ est nul. Nous avons donc : 
\begin{eqnarray}
L_kL_k\overline{L}_k f_n&=&\sum_{j=1}^{n-1}\sum_{l=1}^{n-1} \overline{b}_j^l \left(L_kL_kf_l\overline{L}_k \overline{f}_j-\gamma_{k,\overline{k}}L_kf_lT \overline{f}_j-\gamma_{k,\overline{k}} L_kf_lT \overline{f}_j\right)\nonumber\\
&=&\sum_{j=1}^{n-1}\sum_{l=1}^{n-1} \overline{b}_j^l \left(L_kL_kf_l\overline{L}_k \overline{f}_j-2\gamma_{k,\overline{k}}L_kf_lT \overline{f}_j\right).\label{Lkkkbfn}\\
&\in& \overline{\mathcal{C}}_{2,1}\nonumber
\end{eqnarray}
Ainsi, $TL_kf_n =\dfrac{1}{2\gamma_{\overline{k},k}}\left(L_k L_k\overline{L}_k  f_n-\overline{L}_kL_kL_k f_n\right)\in \overline{\mathcal{C}}_{2,1}$.
\end{itemize}

\underline{\textit{(iii) Montrons que $T L_kL_m f_p \in \overline{\mathcal{C}}_{2,1}$}}\\
Nous avons démontré en (\ref{l1.3}) que $L_kL_kL_m f_p \in \overline{\mathcal{C}}_{1,1}$. Grâce au Lemme \ref{l00}, nous avons $\overline{L}_kL_kL_kL_m f_p \in \overline{\mathcal{C}}_{2,1}$.\\
Pour calculer le terme $TL_kL_m f_p$, nous partons de $\overline{L}_kL_kL_kL_mf_p $ et nous utilisons trois fois l'égalité (\ref{rho}).
\begin{eqnarray*}
\overline{L}_kL_kL_kL_m f_p&=& L_k\overline{L}_k L_kL_m f_p -  [L_k,\overline{L}_k] L_kL_mf_p \\ 
&=&L_kL_k\overline{L}_k L_m f_p - L_k[L_k,\overline{L}_k] L_m f_p - \gamma_{k,\overline{k}}T L_kL_m f_p \\
&=&L_kL_k L_m \overline{L}_k f_p-L_kL_k(\gamma_{m,\overline{k}}T)  f_p- L_k(\gamma_{k,\overline{k}}T) L_m f_p - \gamma_{k,\overline{k}}T L_kL_m f_p  
\end{eqnarray*}
Les termes $\gamma_{m,\overline{k}}$ et $\gamma_{k,\overline{k}}$ étant constants, il vient : 
\begin{eqnarray*}
\overline{L}_kL_kL_kL_m f_p&=&L_kL_k L_m \overline{L}_k f_p- \gamma_{m,\overline{k}}L_kL_kT  f_p-\gamma_{k,\overline{k}}L_kT L_m f_p - \gamma_{k,\overline{k}}T L_kL_m f_p  
\end{eqnarray*}
Les champs $T$ et $L_{k}$ commutant, nous obtenons : 
\begin{eqnarray*}
\overline{L}_kL_kL_kL_m f_p&=&L_kL_k L_m \overline{L}_k f_p- \gamma_{m,\overline{k}}TL_kL_k  f_p- 2\gamma_{k,\overline{k}}T L_kL_m f_p.  
\end{eqnarray*}

\begin{itemize}
\item \underline{Si $p\neq n$ :} \\
D'après (\ref{3.1}), nous avons $\overline{L}_k f_p=0$, et nous obtenons : 
\begin{eqnarray}
\overline{L}_kL_kL_kL_m f_p&=&- \gamma_{m,\overline{k}}TL_kL_k  f_p- 2\gamma_{k,\overline{k}}T L_kL_m f_p.  \label{TLalpha}
\end{eqnarray}
- \underline{Si $m=k$ :} 
Nous avons : $T L_kL_k f_p =\dfrac{-\overline{L}_kL_kL_kL_k f_p}{3\gamma_{k,\overline{k}}}\in \overline{\mathcal{C}}_{2,1}$, pour $p=1,\ldots , n-1$.
- \underline{Si $m\neq k$ :} 
Nous venons de démontrer que $T L_kL_k f_p \in \overline{\mathcal{C}}_{2,1}$. Ainsi, \\
$T L_kL_m f_p =\dfrac{-\overline{L}_kL_kL_kL_m f_p- \gamma_{m,\overline{k}}TL_kL_k  f_p}{2\gamma_{k,\overline{k}}}\in \overline{\mathcal{C}}_{2,1}$, pour $p=1,\ldots , n-1$.

\item \underline{Si $p=n$.}\\
Nous avons : 
\begin{equation*}\overline{L}_kL_kL_kL_m f_p = L_kL_k L_m \overline{L}_k f_n - \gamma_{m,\overline{k}}TL_kL_k  f_n - 2\gamma_{k,\overline{k}}T L_kL_m f_p.\\ \end{equation*}
Calculons le terme $ L_kL_kL_m \overline{L}_k f_n$. Appliquons les champs $L_kL_kL_m$ à l'égalité (\ref{Lkbfn}). Nous obtenons :
\begin{eqnarray*}
 L_kL_kL_m \overline{L}_k f_n&=& L_kL_kL_m\left(\sum_{j=1}^{n-1}\overline{\beta}_j(f)\overline{L}_k \overline{f}_j\right)
\end{eqnarray*}
Remplaçons $\overline{\beta}_j$ par sa valeur donnée dans (\ref{betai}) :
\begin{eqnarray*}
 L_kL_kL_m \overline{L}_k f_n&=& L_kL_kL_m\left(\sum_{j=1}^{n-1}(\sum_{l=1}^{n-1}\overline{a}_j^l \overline{f}_l+\overline{b}_j^l f_l)\overline{L}_k \overline{f}_j\right)\\
&=& L_kL_k\left(\sum_{j=1}^{n-1}\sum_{l=1}^{n-1}\overline{b}_j^l\left( L_mf_l\overline{L}_k \overline{f}_j+f_lL_m\overline{L}_k \overline{f}_j\right)+\overline{a}_j^l \overline{f}_lL_m\overline{L}_k \overline{f}_j\right).\\
\end{eqnarray*}
Or, $L_m\overline{L}_k \overline{f}_j=\gamma_{m,\overline{k}}T\overline{f}_j$. Nous avons donc :  
\begin{eqnarray*}
 L_kL_kL_m \overline{L}_k f_n&=& L_kL_k\left(\sum_{j=1}^{n-1}\sum_{l=1}^{n-1}\overline{b}_j^l\left( L_mf_l\overline{L}_k \overline{f}_j+\gamma_{m,\overline{k}}f_lT\overline{f}_j\right)+\gamma_{m,\overline{k}}\overline{a}_j^l \overline{f}_lT \overline{f}_j\right)\\
&=& L_k\left(\sum_{j=1}^{n-1}\sum_{l=1}^{n-1}\overline{b}_j^l\left( L_kL_mf_l\overline{L}_k \overline{f}_j +L_mf_lL_k\overline{L}_k \overline{f}_j+\gamma_{m,\overline{k}}\left(L_mf_lT \overline{f}_j+f_lL_kT\overline{f}_j\right)+0\right)\right). \\
\end{eqnarray*}
Or, $L_k\overline{L}_k \overline{f}_j=-\gamma_{\overline{k},k}T\overline{f}_j$. Ainsi, 
\begin{eqnarray*}
 L_kL_kL_m \overline{L}_k f_n&=& L_k\left(\sum_{j=1}^{n-1}\sum_{l=1}^{n-1}\overline{b}_j^l\left( L_kL_mf_l\overline{L}_k \overline{f}_j +\left(\gamma_{m,\overline{k}}-\gamma_{\overline{k},k}\right)L_mf_lT \overline{f}_j\right)\right) \\
&=& \sum_{j=1}^{n-1}\sum_{l=1}^{n-1}\overline{b}_j^l\left(L_k L_kL_mf_l\overline{L}_k \overline{f}_j +L_kL_mf_lL_k\overline{L}_k \overline{f}_j\right.\\
&&\left.\hspace{1,4cm}+\left(\gamma_{m,\overline{k}}-\gamma_{\overline{k},k}\right)\left(L_kL_mf_lT \overline{f}_j+L_mf_lL_kT \overline{f}_j\right)\right).
\end{eqnarray*}
Or, les champs $T$ et $L_k$ commutent et $L_k\overline{L}_k \overline{f}_j=-\gamma_{\overline{k},k}T\overline{f}_j$. Ainsi,
\begin{eqnarray*}
 L_kL_kL_m \overline{L}_k f_n&=&\sum_{j=1}^{n-1}\sum_{l=1}^{n-1}\overline{b}_j^l\left(L_k L_kL_mf_l\overline{L}_k \overline{f}_j -\gamma_{\overline{k},k}L_kL_mf_lT \overline{f}_j+\left(\gamma_{m,\overline{k}}-\gamma_{\overline{k},k}\right)L_kL_mf_lT \overline{f}_j\right)\\
&=&\sum_{j=1}^{n-1}\sum_{l=1}^{n-1}\overline{b}_j^l\left(L_k L_kL_mf_l\overline{L}_k \overline{f}_j +\left(\gamma_{m,\overline{k}}-2\gamma_{\overline{k},k}\right)L_kL_mf_lT \overline{f}_j\right)\\
&\in&\overline{\mathcal{C}}_{2,1}.\\
\end{eqnarray*}
- \underline{Si $m=k$ :} Reprenons le calcul en (\ref{TLalpha}) :
\begin{eqnarray*}
 T L_kL_k f_n=\dfrac{1}{3\gamma_{k,\overline{k}}}( L_kL_k L_m \overline{L}_k f_n-\overline{L}_kL_kL_kL_m f_n )&\in&\overline{\mathcal{C}}_{2,1}.\\
\end{eqnarray*}
- \underline{Si $m\neq k$ :} Nous avons, toujours d'après l'égalité (\ref{TLalpha})  :

\begin{eqnarray*}
 T L_kL_k f_n=\dfrac{1}{ 2\gamma_{k,\overline{k}}}( L_kL_k L_m \overline{L}_k f_n-\overline{L}_kL_kL_kL_m f_n+\gamma_{m,\overline{k}}TL_kL_k  f_n  )&\in&\overline{\mathcal{C}}_{2,1}.\\
\end{eqnarray*}
\end{itemize}

\textit{(iv) } Nous pouvons démontrer que $T L_kL_m L_jf_p \in \overline{\mathcal{C}}_{2,1}$ en utilisant la même technique que précédemment.\\
Ceci achève la démonstration de la Proposition \ref{l2}.$\square$
\begin{rem}
A ce stade, nous avons démontré l'assertion $(ii)$ de la Proposition \ref{l4}. Il reste à traiter les cas $t=2$ et $t=3$ pour l'assertion $(i)$.
\end{rem}
\hspace{1cm}\

Nous pouvons maintenant démontrer la Proposition \ref{l3}, qui constitue la troisième étape de la démonstration de la Proposition \ref{l4} :
\begin{prop}\label{l3} Pour $p=1,\ldots ,n $, pour $k=1,\ldots ,n-1 $ : 
\begin{eqnarray*}
(i)&T² f_p \in \overline{\mathcal{C}}_{3,1}, \\
(ii)&T² L_k f_p \in \overline{\mathcal{C}}_{3,1}.
\end{eqnarray*}
\end{prop}
\noindent\textsc{Preuve de la Proposition \ref{l3}} : Le calcul est similaire à celui du (iii) de la Proposition \ref{l2} : \\
\underline{\textit{ (i) Montrons que $T² f_p \in \overline{\mathcal{C}}_{3,1}$ :}}\\
Nous avons déjà démontré que $T L_1 f_p \in \overline{\mathcal{C}}_{2,1}$. D'après le Lemme \ref{l00}, nous obtenons que $\overline{L}_1T L_1  f_p \in \overline{\mathcal{C}}_{3,1}$. \\
Les champs $T$ et $L_1$ commutant, nous pouvons écrire : 
\begin{eqnarray*}\overline{L}_1T L_1  f_p&=&T \overline{L}_1 L_1 f_p. \end{eqnarray*}
Nous avons déjà calculé $\overline{L}_1 L_1 f_p$ à l'égalité (\ref{L1bL1fp}) :
\begin{eqnarray*}
\overline{L_1}T L_1  f_p& = & T ( L_1\overline{L_1} f_p -\gamma_{\overline{1},1}Tf_p) \\
& = & T  L_1 \overline{L_1}f_p- \gamma_{1, \overline{1}}T²f_p.
\end{eqnarray*}
\begin{itemize}
 \item \underline{Si $p\neq n$}: \\
Nous avons  $\overline{L}_1  f_p=0$. Ainsi, 
\begin{eqnarray*}
\overline{L}_1T L_1 f_p&=&-\gamma_{1,\overline{1}}T²f_p.
\end{eqnarray*}
Ainsi, $T²f_p=\dfrac{-1}{\gamma_{1,\overline{1}}}\overline{L}_1T L_1 f_p\in \overline{\mathcal{C}}_{3,1}$.

\item \underline{Si $p= n$}: \\
Nous avons
\begin{eqnarray*}
\overline{L}_1T L_1 f_p&=&TL_1\overline{L}_1   f_n+T(-\gamma_{1,\overline{1}}Tf_n)\\
&=&TL_1\overline{L}_1   f_n-\gamma_{1,\overline{1}}T²f_n.
\end{eqnarray*}
Or, en remplaçant $L_1\overline{L}_1   f_n$ par sa valeur donnée dans (\ref{L11fn}), nous obtenons : 
\begin{eqnarray}
    TL_1\overline{L}_1   f_n & = & T\left(\sum_{j=1}^{n-1}\left(\sum_{l=1}^{n-1}\overline{b}_j^l \left(L_1f_l\overline{L}_1 \overline{f}_j-\gamma_{1,\overline{1}}f_lT \overline{f}_j\right)-\gamma_{1,\overline{1}}\overline{a}_j^l\overline{f}_lT \overline{f}_j\right)\right)\nonumber\\
 & = & \sum_{j=1}^{n-1}\sum_{l=1}^{n-1}\left(\overline{b}_j^l \left(TL_1f_l\overline{L}_1 \overline{f}_j+L_1f_lT\overline{L}_1 \overline{f}_j-\gamma_{1,\overline{1}}\left(Tf_lT \overline{f}_j+f_lT^2 \overline{f}_j\right)\right)\right.\nonumber\\
&&\left.\hspace{1,4cm}-\gamma_{1,\overline{1}}\overline{a}_j^l\left(T\overline{f}_lT \overline{f}_j+\overline{f}_lT^2 \overline{f}_j\right)\right).\label{Tl11bfn}
   \end{eqnarray}
Tous les termes qui apparaissent dans le membre de droite de l'égalité (\ref{Tl11bfn}) sont dans l'espace $\overline{\mathcal{C}}_{2,1}$, sauf $T^2\overline{f}_j$. Or nous venons de démontrer que le terme $T²f_j $ appartient à l'espace  $\overline{\mathcal{C}}_{3,1}$. Par conjugaison, nous avons : $T²\overline{f}_j \in \mathcal{C}_{3,1}$. \\
Le terme $T²\overline{f}_j $ s'écrit donc comme une fonction analytique réelle de $\{\,T^tL^\alpha f_j : t+|\alpha|\leqslant 3\, , t\leqslant 1 ,\, j=1,\ldots,n-1,\; L^\alpha f_n : |\alpha|\leqslant 3\}.$ Mais nous avons démontré dans les Propositions \ref{l1} et \ref{l2} que ces termes sont dans l'espace $\overline{\mathcal{C}}_{2,1}$. Nous pouvons donc conclure que le terme $T²\overline{f}_j $ appartient à l'espace $  \overline{\mathcal{C}}_{2,1}$.\\
Ainsi, $T²f_n=\dfrac{-1}{\gamma_{1,\overline{1}}}(\overline{L}_1T L_1 f_n-TL_1\overline{L}_1   f_n)\in \overline{\mathcal{C}}_{3,1}$.\\
\end{itemize}

\underline{\textit{(ii) Montrons  que $T² L_k f_p \in \overline{\mathcal{C}}_{3,1}$ : }}\\
Nous avons démontré que $T L_k L_k f_p \in \overline{\mathcal{C}}_{2,1}$. D'après le Lemme \ref{l00}, nous obtenons que $\overline{L}_kT L_k L_k  f_p \in \overline{\mathcal{C}}_{3,1}$. \\
Les champs $T$ et $L_k$ commutant, nous pouvons écrire : 
\begin{eqnarray*}\overline{L}_kT L_k L_k  f_p&=&T \overline{L}_k L_kL_k f_p. \end{eqnarray*}
Or, nous avons déjà calculé $\overline{L}_k L_kL_k f_p$ en (\ref{TL1}) :
\begin{eqnarray*}
\overline{L}_kT L_k L_k  f_p&=&T (L_k L_k\overline{L}_k  f_p - 2\gamma_{\overline{k},k}TL_k f_p) \\
&=&T  L_kL_k \overline{L}_kf_p-2\gamma_{k,\overline{k}}T²  L_kf_p. \\
\end{eqnarray*}
\begin{itemize}
 \item \underline{Si $p\neq n $ :}\\
Nous avons  $\overline{L}_k f_p=0$.
Ainsi, $T²L_kf_p=\dfrac{-1}{2\gamma_{k,\overline{k}}}\overline{L}_kT L_k L_kf_p\in \overline{\mathcal{C}}_{3,1}$.

\item \underline{Si $p= n$}: Appliquons le champ $T$ à l'égalité (\ref{Lkkkbfn}).
\begin{eqnarray*}
T  L_kL_k \overline{L}_kf_n&=&T\left(\sum_{j=1}^{n-1}\sum_{l=1}^{n-1} \overline{b}_j^l \left(L_kL_kf_l\overline{L}_k \overline{f}_j-2\gamma_{k,\overline{k}}L_kf_lT \overline{f}_j\right)\right)\\
&=&\sum_{j=1}^{n-1}\sum_{l=1}^{n-1} \overline{b}_j^l \left(TL_kL_kf_l\overline{L}_k \overline{f}_j+L_kL_kf_lT\overline{L}_k \overline{f}_j-2\gamma_{k,\overline{k}}\left(TL_kf_lT \overline{f}_j+L_kf_lT² \overline{f}_j\right)\right).
\end{eqnarray*}
Tous les termes qui apparaissent dans le membre de droite sont dans l'espace $\overline{\mathcal{C}}_{2,1}$.
Ainsi, $T²L_kf_n=\dfrac{-1}{2\gamma_{k,\overline{k}}}(\overline{L}_kT L_k L_k f_n-TL_kL_k\overline{L}_k   f_n)\in \overline{\mathcal{C}}_{3,1}$.
\end{itemize}

Nous obtenons ainsi que $T^t L^\alpha f_p \in \overline{\mathcal{C}}_{1+t,1}$ pour $p=1,\ldots, n$, $t+|\alpha|\leqslant 3$ et $t\leqslant 2$, ce qui termine la démonstration de la Proposition \ref{l3}.$\square$\

Pour compléter la démonstration de la Proposition \ref{l1}, il reste à démontrer que :
\begin{prop}\label{l5} Pour $p=1,\ldots ,n $ : 
\begin{eqnarray*}
T³ f_p \in \overline{\mathcal{C}}_{4,1} .
\end{eqnarray*}
\end{prop}
\noindent\textsc{Preuve de la Proposition \ref{l5}} : Nous avons déjà démontré que $T² L_1 f_p \in \overline{\mathcal{C}}_{3,1}$. D'après le Lemme \ref{l00}, nous obtenons que $\overline{L}_1T^{2} L_1  f_p \in \overline{\mathcal{C}}_{4,1}$. Les champs $T$ et $L_1$ commutant, nous pouvons écrire : 
\begin{eqnarray*}\overline{L}_1T^{2} L_1 f_p&=&T² \overline{L}_1 L_1  f_p. \end{eqnarray*}
Nous avons déjà calculé $\overline{L}_1 L_1  f_p$ en (\ref{L1bL1fp}) : 
\begin{eqnarray*}
\overline{L}_1T^{2} L_1 f_p&=&T²  (L_1\overline{L}_1 f_p -\gamma_{\overline{1},1}Tf_p)\\
&=&T²  L_1\overline{L}_1  f_p-\gamma_{\overline{1},1}T³ f_p. \\
\end{eqnarray*}
\begin{itemize}
 \item \underline{Si $p\neq n$}: \\
Nous avons  $\overline{L}_1 f_p=0$.
Ainsi, $T³f_p=\dfrac{-1}{\gamma_{1,\overline{1}}}\overline{L}_1T² L_1 f_p\in \overline{\mathcal{C}}_{4,1}$.

\item \underline{Si $p= n$}: \\
Appliquons le champ $T$ dans l'égalité (\ref{Tl11bfn}) :
\begin{eqnarray}
    T²L_1\overline{L}_1   f_n&=&T\left(\sum_{j=1}^{n-1}\sum_{l=1}^{n-1}\overline{b}_j^l \left(TL_1f_l\overline{L}_1 \overline{f}_j+L_1f_lT\overline{L}_1 \overline{f}_j-\gamma_{1,\overline{1}}\left(Tf_lT \overline{f}_j+f_lT² \overline{f}_j\right)\right)\right.\nonumber\\
&&\hspace{1,8cm}\left.-\gamma_{1,\overline{1}}\overline{a}_j^l\left(T\overline{f}_lT \overline{f}_j+\overline{f}_lT^2 \overline{f}_j\right)\right)\nonumber\\
&=&\sum_{j=1}^{n-1}\sum_{l=1}^{n-1}\left(\overline{b}_j^l \left(T²L_1f_l\overline{L}_1 \overline{f}_j+2TL_1f_lT\overline{L}_1 \overline{f}_j+L_1f_lT²\overline{L}_1 \overline{f}_j\right.\right.\nonumber\\
 &&\hspace{1,3cm}\left.\left.-\gamma_{1,\overline{1}}\left(T²f_lT \overline{f}_j+Tf_lT² \overline{f}_j+f_lT³ \overline{f}_j\right)\right)\right.\nonumber\\
&&\left.\hspace{1,3cm}-\gamma_{1,\overline{1}}\overline{a}_j^l\left(T^2\overline{f}_lT \overline{f}_j+2T\overline{f}_lT^2 \overline{f}_j+\overline{f}_lT^3 \overline{f}_j\right)\right).\label{T2L1L1bfn}
\end{eqnarray}

Tous les termes qui apparaissent dans le membre de droite de l'égalité (\ref{T2L1L1bfn}) sont dans l'espace $\overline{\mathcal{C}}_{3,1}$, sauf $T^3\overline{f}_j$. Or nous venons de démontrer que le terme que $T³f_j$ appartient à l'espace $\overline{\mathcal{C}}_{4,1}$. Par conjugaison, nous avons :  $T³\overline{f}_j \in \mathcal{C}_{4,1}$.\\
Le terme $T³\overline{f}_j $ s'écrit donc comme une fonction analytique réelle de $\{\,T^tL^\alpha f_j : t+|\alpha|\leqslant 4\, , t\leqslant 1 ,\, j=1,\ldots,n-1,\; L^\alpha f_n : |\alpha|\leqslant 4\}.$ Mais nous avons démontré dans les Propositions \ref{l1} et \ref{l2} que ces termes sont dans l'espace $\overline{\mathcal{C}}_{2,1}$. Nous pouvons donc conclure donc que le terme $T³\overline{f}_j$ appartient à l'espace $  \overline{\mathcal{C}}_{2,1}$.\\
Ainsi, $T³f_n=\dfrac{-1}{\gamma_{1,\overline{1}}}(\overline{L}_1T² L_1 f_n-T²L_1\overline{L}_1   f_n)\in \overline{\mathcal{C}}_{4,1}$.
\end{itemize}
$\square$\\
Ceci termine la preuve de la Proposition \ref{l4} $\square$.

\subsection{Démonstration du Théorème \ref{prop1}.}
Pour obtenir un système complet vérifié par la fonction $f$, et ainsi terminer la démonstration du Théorème \ref{prop1}, nous procédons en trois étapes. Nous commençons d'abord, dans le Lemme \ref{l2.3.1}, par exprimer les dérivées d'ordre $3$ de $f$ de la forme $T^tL^\alpha f_j$ comme des fonctions analytiques réelles des dérivées d'ordre inférieur ou égal à 2. Nous utilisons pour cela le résultat obtenu dans la Proposition \ref{l4}. Puis, dans le Lemme \ref{beta}, nous exprimons les dérivées de la forme $T^tL^\alpha \overline{L}^\beta f_j$. Enfin, nous exprimons toutes les dérivées dans la Proposition \ref{crochet}.\\
Démontrons maintenant le lemme correspondant à la première étape : 
\begin{lem}\label{l2.3.1}Pour $t+|\alpha|=3$, pour $j=1,\ldots,n-1$,
\begin{equation*}
T^tL^\alpha f_j \in \mathcal{C}_{2,1}.
\end{equation*}
\end{lem}\

\noindent\textsc{Preuve du Lemme \ref{l2.3.1}} : La preuve de ce lemme est encore basée sur \cite{H97}. La Proposition \ref{l4} signifie que pour $t+|\alpha|\leqslant p$, pour $t\leqslant q$, pour $j=1,\ldots,n-1$,
\begin{equation}\label{3.15}
T^tL^\alpha f_j \in \overline{\mathcal{C}}_{1+q,1},\text{ si } q\leqslant p\leqslant 3.
\end{equation}

Or, l'assertion $(ii)$ de la Proposition \ref{l4} implique l'inclusion suivante : 
\begin{equation}\label{3.17}
\mathcal{C}_{4,1}\subset \overline{\mathcal{C}_{2,1}}.
\end{equation}

En prenant le conjugué (\ref{3.17}), nous obtenons : 
\begin{equation}\label{3.18}
\overline{\mathcal{C}_{4,1}}\subset \mathcal{C}_{2,1}.
\end{equation}

En prenant $p=q=3$ dans (\ref{3.15}), et en utilisant l'inclusion (\ref{3.18}), nous avons : 
\begin{equation*}
T^tL^\alpha f_j \in \mathcal{C}_{2,1}.
\end{equation*}\\
Ceci termine la démonstration du Lemme \ref{l2.3.1}.$\square$\\
En particulier, nous avons aussi
\begin{equation*}
\mathcal{C}_{3,3} \subset \mathcal{C}_{2,1}.
\end{equation*}

Nous pouvons maintenant obtenir les dérivées d'ordre 3 de la forme $T^tL^\alpha \overline{L}^\beta f_j$ comme des fonctions analytiques réelles des dérivées d'ordre inférieur ou égal à deux : 
\begin{lem} \label{beta}Les dérivées de la forme $T^tL^\alpha \overline{L}^\beta f_j$, (où $t+|\alpha|+|\beta|=3$) appartiennent à l'espace $\mathcal{C}_{2,1}$ et s'expriment donc comme une fonction analytique réelle des dérivées d'ordre inférieur ou égal à 2. \\
\end{lem}
\noindent\textsc{Preuve du Lemme \ref{beta} : }Les deux éléments importants de la preuve de ce lemme sont le résultat du Lemme \ref{l2.3.1} ainsi que la relation (\ref{3.1bis}).
\begin{itemize}
\item Si $|\beta|=0$, d'après la Proposition \ref{l2.3.1}, nous savons que $T^tL^\alpha f_j$ est une fonction analytique réelle de $\{T^sL^\beta f, \; s+|\beta|\leqslant 2,\;s\leqslant 1\}$.
\item Si $|\beta|\neq0$ et $j \neq n$, alors $T^tL^\alpha \overline{L}^\beta f_j=0$ d'après (\ref{3.1}).\\
\item Si $|\beta|\neq0$ et $j = n$, alors nous avons six formes de termes possibles : $T^2 \overline{L}_k f_n$,  $TL_m \overline{L}_k f_n$, $T\overline{L}_m \overline{L}_k f_n$, $L_p L_m \overline{L}_k f_n$, $L_p\overline{L}_m \overline{L}_k f_n$ et $\overline{L}_p\overline{L}_m \overline{L}_k f_n$.
\begin{itemize}
\item En appliquant les champs $T²$ à l'égalité (\ref{Lkbfn}), nous obtenons :
\begin{eqnarray*}
 T²\overline{L}_k f_n&=&T\sum_{j=1}^{n-1}\sum_{l=1}^{n-1}\left(\left(\overline{b}_j^l Tf_l+\overline{a}_j^l T\overline{f}_l\right)\overline{L}_k \overline{f}_j+\left(\overline{b}_j^l f_l+\overline{a}_j^l \overline{f}_l\right)T\overline{L}_k \overline{f}_j\right)\\
&=&\sum_{j=1}^{n-1}\sum_{l=1}^{n-1}\left(\left(\overline{b}_j^l T^2f_l+\overline{a}_j^l T^2\overline{f}_l\right)\overline{L}_k \overline{f}_j+2\left(\overline{b}_j^l Tf_l+\overline{a}_j^l T\overline{f}_l\right)T\overline{L}_k \overline{f}_j\right.\\
&&\hspace{2cm}\left.+\left(\overline{b}_j^l f_l+\overline{a}_j^l \overline{f}_l\right)T^2\overline{L}_k \overline{f}_j\right)\\
&\in&\overline{\mathcal{C}}_{3,2}\subset \mathcal{C}_{3,1}\subset\mathcal{C}_{2,1}\text{ d'après le conjugué de (\ref{3.15}) et (\ref{3.17}).}
\end{eqnarray*}
Nous avons donc $T^2\overline{L}_k f_n \in \mathcal{C}_{2,1}$.\\

\item En appliquant le champ $L_m$ à l'égalité (\ref{Lkbfn}), nous obtenons :\\
 \begin{eqnarray}L_m\overline{L}_k f_n&=&\sum_{j=1}^{n-1}\sum_{l=1}^{n-1}
\left(\overline{b}_j^l L_mf_l\overline{L}_k \overline{f}_j+\left(\overline{b}_j^l f_l+\overline{a}_j^l \overline{f}_l\right)L_m\overline{L}_k \overline{f}_j\right)\notag\\
  &=&\sum_{j=1}^{n-1}\sum_{l=1}^{n-1}\left(\overline{b}_j^l L_mf_l\overline{L}_k \overline{f}_j-\gamma_{\overline{k},m}\left(\overline{b}_j^l f_l+\overline{a}_j^l \overline{f}_l\right)T \overline{f}_j\right).
\label{LmLkbfn}
\end{eqnarray}
Donc, en appliquant le champ $T$ à l'égalité précédente : 
\begin{eqnarray*}TL_m\overline{L}_k f_n&=&\sum_{j=1}^{n-1}\sum_{l=1}^{n-1}\left(\overline{b}_j^l \left(TL_mf_l\overline{L}_k \overline{f}_j+L_mf_lT\overline{L}_k \overline{f}_j\right)\right.\\
&&\hspace{1cm}\left.-\gamma_{\overline{k},m}\left(\left(\overline{b}_j^l Tf_l+\overline{a}_j^l T\overline{f}_l\right)T \overline{f}_j+\left(\overline{b}_j^l f_l+\overline{a}_j^l \overline{f}_l\right)T^2 \overline{f}_j\right)\right)\\
&&\in \overline{\mathcal{C}_{2,2}} 
\subset\mathcal{C}_{2,1}.\end{eqnarray*}
Enfin, en appliquant le champ $L_p$ à l'égalité (\ref{LmLkbfn}), nous obtenons : 
\begin{eqnarray*} L_pL_m\overline{L}_k f_n&=&\sum_{j=1}^{n-1}\sum_{l=1}^{n-1}
\left(\overline{b}_j^l (L_pL_mf_l\overline{L}_k \overline{f}_j+L_mf_lL_p\overline{L}_k \overline{f}_j)-\gamma_{\overline{k},m}\overline{b}_j^l L_pf_lT \overline{f}_j\right)\\
&=&\sum_{j=1}^{n-1}\sum_{l=1}^{n-1}
\left(\overline{b}_j^l (L_pL_mf_l\overline{L}_k \overline{f}_j-\gamma_{\overline{k},p}L_mf_lT \overline{f}_j)-\gamma_{\overline{k},m}\overline{b}_j^l L_pf_lT \overline{f}_j\right)\\     
&&\in\mathcal{C}_{2,1}.                                            
 \end{eqnarray*}

\item En appliquant le champ $\overline{L}_m$ à l'égalité (\ref{Lkbfn}), nous obtenons : 
\begin{eqnarray}
\overline{L}_m\overline{L}_k f_n&=&\sum_{j=1}^{n-1}\sum_{l=1}^{n-1}\left(\overline{a}_j^l \overline{L}_m\overline{f}_l\overline{L}_k \overline{f}_j+\left(\overline{b}_j^l f_l+\overline{a}_j^l \overline{f}_l\right)\overline{L}_m\overline{L}_k \overline{f}_j\right).
\label{LmbLkbfn}
\end{eqnarray}
Donc, en appliquant le champ $T$ à l'égalité précédente :
\begin{eqnarray*}T\overline{L}_m\overline{L}_k f_n&=&\sum_{j=1}^{n-1}\sum_{l=1}^{n-1}\left(\overline{a}_j^l\left(T \overline{L}_m\overline{f}_l\overline{L}_k \overline{f}_j+\overline{L}_m\overline{f}_lT\overline{L}_k \overline{f}_j\right)+\left(\overline{b}_j^l Tf_l+\overline{a}_j^l T\overline{f}_l\right)\overline{L}_m\overline{L}_k \overline{f}_j\right.\\
&&\hspace{1,5cm}\left.+\left(\overline{b}_j^l f_l+\overline{a}_j^l \overline{f}_l\right)T\overline{L}_m\overline{L}_k \overline{f}_j\right)\\
&&\in \overline{\mathcal{C}_{3,1}} \subset\mathcal{C}_{2,1}.
\end{eqnarray*}

De plus, appliquons le champ $L_p$ à l'égalité (\ref{LmbLkbfn}) : 
\begin{eqnarray*}
L_p\overline{L}_m\overline{L}_k f_n&=&\sum_{j=1}^{n-1}\sum_{l=1}^{n-1}
\left(\overline{a}_j^l \left(L_p\overline{L}_m\overline{f}_l\overline{L}_k \overline{f}_j+\overline{L}_m\overline{f}_lL_p\overline{L}_k \overline{f}_j\right)+\overline{b}_j^l L_pf_l\overline{L}_m\overline{L}_k \overline{f}_j\right.\\
&&\hspace{1,5cm}\left.+\left(\overline{b}_j^l f_l+\overline{a}_j^l \overline{f}_l\right)L_p\overline{L}_m\overline{L}_k \overline{f}_j\right)\\
&=&\sum_{j=1}^{n-1}\sum_{l=1}^{n-1}
\left(\overline{a}_j^l \left(\gamma_{p,\overline{m}}T\overline{f}_l\overline{L}_k \overline{f}_j+\gamma_{p,\overline{k}}\overline{L}_m\overline{f}_lT \overline{f}_j\right)+\overline{b}_j^l L_pf_l\overline{L}_m\overline{L}_k \overline{f}_j\right.\\
&&\hspace{1,5cm}\left.+\left(\overline{b}_j^l f_l+\overline{a}_j^l \overline{f}_l\right)\left(\overline{L}_mL_p\overline{L}_k \overline{f}_j+\gamma_{p,\overline{m}}T\overline{L}_k \overline{f}_j\right)\right)\\
&=&\sum_{j=1}^{n-1}\sum_{l=1}^{n-1}
\left(\overline{a}_j^l \left(\gamma_{p,\overline{m}}T\overline{f}_l\overline{L}_k \overline{f}_j+\gamma_{p,\overline{k}}\overline{L}_m\overline{f}_lT \overline{f}_j\right)+\overline{b}_j^l L_pf_l\overline{L}_m\overline{L}_k \overline{f}_j\right.\\
&&\hspace{1,5cm}\left.+\left(\overline{b}_j^l f_l+\overline{a}_j^l \overline{f}_l\right)\left(\gamma_{p,\overline{k}}\overline{L}_mT \overline{f}_j+\gamma_{p,\overline{m}}T\overline{L}_k \overline{f}_j\right)\right)\\
&=&\sum_{j=1}^{n-1}\sum_{l=1}^{n-1}
\left(\overline{a}_j^l \left(\gamma_{p,\overline{m}}T\overline{f}_l\overline{L}_k \overline{f}_j+\gamma_{p,\overline{k}}\overline{L}_m\overline{f}_lT \overline{f}_j\right)+\overline{b}_j^l L_pf_l\overline{L}_m\overline{L}_k \overline{f}_j\right.\\
&&\hspace{1,5cm}\left.+\left(\overline{b}_j^l f_l+\overline{a}_j^l \overline{f}_l\right)\left(\gamma_{p,\overline{k}}T\overline{L}_m \overline{f}_j+\gamma_{p,\overline{m}}T\overline{L}_k \overline{f}_j\right)\right)\\
&&\in \overline{\mathcal{C}_{2,1}} \subset \mathcal{C}_{2,1}.
\end{eqnarray*}
Enfin, appliquons le champ $\overline{L}_p$ à l'égalité (\ref{LmbLkbfn}) :
\begin{eqnarray*}
\overline{L}_p\overline{L}_m\overline{L}_k f_n &=& \sum_{j=1}^{n-1}\sum_{l=1}^{n-1}\left(\overline{a}_j^l\left( \overline{L}_p\overline{L}_m\overline{f}_l\overline{L}_k \overline{f}_j+\overline{L}_m\overline{f}_l\overline{L}_p\overline{L}_k \overline{f}_j\right)+\overline{a}_j^l\overline{L}_p \overline{f}_l\overline{L}_m\overline{L}_k \overline{f}_j\right.\\
&&\hspace{1,5cm}\left.+\left(\overline{b}_j^l f_l+\overline{a}_j^l \overline{f}_l\right)\overline{L}_p\overline{L}_m\overline{L}_k \overline{f}_j\right)\\
&&\in \overline{\mathcal{C}_{3,0}} \subset\mathcal{C}_{2,1}.
\end{eqnarray*}
\end{itemize}
\end{itemize}
Ceci termine la démonstration du Lemme \ref{beta}.
$\square$\\

Enfin, nous pouvons maintenant exprimer toutes les dérivées d'ordre 3 de $f$ comme des fonctions analytiques réelles des dérivées d'ordre inférieur ou égal à deux : 
\begin{prop}\label{crochet}
Toutes les dérivées d'ordre 3 de $f=(f_1,\ldots,f_{n})$ s'expriment comme des fonctions analytiques réelles des dérivées d'ordre inférieur ou égal à 2.
\end{prop}

\noindent\textsc{Preuve de la Proposition \ref{crochet}} : D'après le Lemme \ref{beta}, il s'agit de démontrer que les dérivées d'ordre 3 qui ne sont pas de la forme $T^tL^\alpha\overline{L}^\beta f_j$ s'expriment comme des fonctions analytiques réelles des dérivées d'ordre inférieur ou égal à 2.\\
Pour une telle dérivée, soit $t$ le nombre de fois où le champ de vecteurs $T$ apparaît, $a$ le nombre de fois où un champ $L_j$ apparaît, et $b$ le nombre de fois où un champ $\overline{L}_j$ apparaît. 
Il s'agit essentiellement d'effectuer des crochets de Lie pour se ramener à des dérivées qui ont déjà été exprimées dans le Lemme \ref{beta}.

\begin{itemize}
 \item \underline{Si $t=2$ et $a=1$ ou si $t=1$ et $a=2$}, les champs $T$ et $L_k$ commutant, la dérivée est égale à une dérivée de la forme $T^2L_kf_j$ ou $TL_kL_mf_j$, qui s'expriment comme des fonctions analytiques réelles des dérivées d'ordre $\leqslant$ 2.
 \item \underline{Si $t=2$ et $b=1$ ou si $t=1$ et $b=2$}, les champs $T$ et $\overline{L}_k$ commutant, la dérivée est égale à une dérivée de la forme $T^2\overline{L}_kf_j$ ou $T\overline{L}_k\overline{L}_mf_j$ qui s'expriment comme des fonctions analytiques réelles des dérivées d'ordre $\leqslant$ 2.
 \item \underline{Si $a=2$ et $b=1$} : il y a deux formes de termes possibles : $\overline{L}_kL_mL_pf_j$ et $L_m\overline{L}_kL_pf_j$.\\
Or, d'après l'égalité (\ref{rho}), nous avons :\\ 
$L_m\overline{L}_kL_pf_j = L_m(L_p\overline{L}_kf_j - [L_p\overline{L}_k]f_j)=L_mL_p\overline{L}_kf_j - \gamma_{p,\overline{k}}L_mTf_j$.\\
Nous pouvons donc exprimer $L_m\overline{L}_kL_pf_j$ comme une fonction analytique réelle des dérivées d'ordre $\leqslant$ 2.\\
Nous avons aussi,\\
$\overline{L}_kL_mL_pf_j= L_m\overline{L}_kL_pf_j-[L_m,\overline{L}_k]L_pf_j=L_m\overline{L}_kL_pf_j-\gamma_{m,\overline{k}}TL_pf_j$.\\
Nous pouvons donc exprimer $\overline{L}_kL_mL_pf_j$ comme une fonction analytique réelle des dérivées d'ordre $\leqslant$ 2.
 \item \underline{Si $a=1$ et $b=2$} : il y a encore deux formes de termes possibles : $L_k\overline{L}_m\overline{L}_pf_j$ et $\overline{L}_mL_k\overline{L}_pf_j$. \\
Or, d'après l'égalité (\ref{rho}), nous avons :\\ 
$\overline{L}_mL_k\overline{L}_pf_j=L_k\overline{L}_m\overline{L}_pf_j-\gamma_{k,\overline{m}}T\overline{L}_pf_j$. \\
Nous pouvons donc exprimer $\overline{L}_mL_k\overline{L}_pf_j$ comme une fonction analytique réelle des dérivées d'ordre $\leqslant$ 2.\\
Nous avons aussi, toujours d'après l'égalité (\ref{rho}),\\
$L_k\overline{L}_m\overline{L}_pf_j=\overline{L}_mL_k\overline{L}_pf_j-\gamma_{\overline{m},k}T\overline{L}_pf_j$.\\
Nous pouvons donc exprimer $L_k\overline{L}_m\overline{L}_pf_j$ comme une fonction analytique réelle des dérivées d'ordre $\leqslant$ 2.
\item \underline{Si $t=1$, $a=1$ et $b=1$} : il y a cinq formes de termes possibles :\\
Les termes $L_mT\overline{L}_kf_j $ et $L_m\overline{L}_kTf_j$ sont tous deux égaux à $TL_m\overline{L}_kf_j$ et s'expriment comme des fonctions analytiques réelles des dérivées d'ordre $\leqslant$ 2.\\ 
De plus, 
$T\overline{L}_kL_mf_j = TL_m\overline{L}_kf_j-\gamma_{m,\overline{k}}T^2f_j$ s'exprime comme une fonction analytique réelle des dérivées d'ordre $\leqslant$ 2 .\\
Enfin, les termes
$\overline{L}_kTL_mf_j $ et $\overline{L}_kL_mTf_j$ sont tous deux égaux à $T\overline{L}_kL_mf_j $. Ce qui termine la preuve de la Proposition \ref{crochet}.$\square$
\end{itemize}

Ainsi, l'application $f$ vérifie un système complet d'ordre 3 avec des applications analytiques réelles. D'après la Proposition \ref{analytique}, l'application $f$ est analytique réelle. Ceci termine la démonstration du Théorème \ref{prop1}.$\square$

\section{Démonstration du Théorème \ref{prop2}}
Nous définissons, dans le premier paragraphe les structures presque complexes vérifiant la condition $(*)$. Dans le deuxième paragraphe, nous démontrons le Théorème \ref{prop2}. Les équations de Cauchy-Riemann tangentielles vérifiées par l'application $g$ sont identiques à celles obtenues dans le cas modèle. Nous pouvons donc conclure la démonstration avec la même méthode que dans le cas modèle. Dans le troisième paragraphe, nous donnons une interprétation géométrique pour les structures presque complexes vérifiant la condition $(*)$.

\subsection{Structures presque complexes vérifiant la condition $(*)$}
\begin{defn} \label{modèle}
Une structure presque complexe $J$ sur $\C^n$ \textbf{vérifie la condition $(*)$} si $J(z)=J_{st}+L(z)$, où $L$ est une matrice $L=(L_{j,k})_{1 \leqslant j,k \leqslant 2n} $ telle que $L_{j,k}=0$ si $1\leqslant j \leqslant 2n-2$, $1\leqslant k \leqslant 2n$.\\
La complexification $J_\C$ d'une structure presque complexe vérifiant la condition $(*)$ s'écrit comme une matrice complexe $2n\times 2n$  :
\begin{eqnarray*}J_{\C}=\begin{pmatrix}
i&0&\ldots&0&0\\
0&-i&\ldots&0&0\\
0&0&\ldots&0&0\\
0&0&\ldots&0&0\\
\ldots&\ldots&\ldots&\ldots&\ldots\\
\tilde{L}_{2n-1,1}(z,\overline{z})&\tilde{L}_{2n-1,2}(z,\overline{z})&\ldots&i+\tilde{L}_{2n-1,2n-1}(z,\overline{z})&\tilde{L}_{2n-1,2n}(z,\overline{z})\\
\overline{\tilde{L}_{2n-1,2}}(z,\overline{z})&\overline{\tilde{L}_{2n-1,1}}(z,\overline{z})&\ldots&\overline{\tilde{L}_{2n-1,2n}}(z,\overline{z})&-i+\overline{\tilde{L}_{2n-1,2n-1}}(z,\overline{z})\\
\end{pmatrix}.\end{eqnarray*}
\end{defn}

La condition $J_{\C}^2=-I_{2n}$ implique que, pour $j=1,\ldots,n-1$, 
\begin{eqnarray}
&& 2i\tilde{L}_{2n-1,2j-1}(z,\overline{z})+\tilde{L}_{2n-1,2j-1}(z,\overline{z})\tilde{L}_{2n-1,2n-1}(z,\overline{z})\notag\\
&&+\overline{\tilde{L}_{2n-1,2j}(z,\overline{z})}\tilde{L}_{2n-1,2n}(z,\overline{z})=0,\label{J1}\\
&&\tilde{L}_{2n-1,2j}(z,\overline{z})\tilde{L}_{2n-1,2n-1}(z,\overline{z})+\overline{\tilde{L}_{2n-1,2j-1}(z,\overline{z})}\tilde{L}_{2n-1,2n}(z,\overline{z})=0.\label{J2}
\end{eqnarray}
Et aussi,
\begin{eqnarray}
 &&2i\tilde{L}_{2n-1,2n-1}(z,\overline{z})+(\tilde{L}_{2n-1,2n-1}(z,\overline{z}))^2+|\tilde{L}_{2n-1,2n}(z,\overline{z})|^2=0,\label{J3}\\
&&\RE(\tilde{L}_{2n-1,2n-1}(z,\overline{z}))=0.\label{J4}
\end{eqnarray}
En particulier, d'après (\ref{J1}), le développement limité à l'ordre $1$ en $0$ de $\tilde{L}_{2n-1,2i-1}$ est nul pour $i=1,\ldots,n-1$.\

Nous étudions maintenant l'espace tangent $J$-holomorphe pour une structure vérifiant la condition $(*)$. 
\begin{lem}\label{tangent}Soit $J$ une structure presque complexe sur $\C^n$ vérifiant la condition $(*)$. Les  champs de vecteurs \begin{eqnarray}
X_i&=&\dzi+a_i(z)\dzn+b_i(z)\dznb,\; i=1,\ldots , n-1,\label{Xi}\\
X_n&=&\dzn+b_n(z)\dznb,\label{Xn}
\end{eqnarray}
 forment une base de l'espace tangent $J$-holomorphe $H^{1,0}_J\C^n$, avec
\begin{eqnarray}
 a_i(z)&=&\dfrac{-i}{2}\tilde{L}_{2n-1,2i-1}(z,\overline{z}), \label{a_i}\\
b_i(z)&=&\dfrac{-i}{2}\overline{\tilde{L}_{2n-1,2i}}(z,\overline{z}),\label{b_i}\\
b_n(z)&=&\dfrac{\overline{\tilde{L}_{2n-1,2n}}(z,\overline{z})}{2i-\overline{\tilde{L}_{2n-1,2n-1}}(z,\overline{z})}.\label{b_n}
\end{eqnarray}
\end{lem}

\noindent\textsc{Preuve du Lemme \ref{tangent}} : On a :

\begin{eqnarray*}
 J X_j&=&J \left(\dzi+a_j(z)\dzn+b_j(z)\dznb\right)\\
&=& i\dzj+\tilde{L}_{2n-1,2j-1}(z,\overline{z})\dzn+\overline{\tilde{L}_{2n-1,2j}(z,\overline{z})}\dznb\\
&&+a_j(z)\left((i+\tilde{L}_{2n-1,2n-1}(z,\overline{z}))\dzn+\overline{\tilde{L}_{2n-1,2n}(z,\overline{z})}\dznb\right)\\
&&+b_j(z)\left(\tilde{L}_{2n-1,2n}(z,\overline{z})\dzn+(-i+\overline{\tilde{L}_{2n-1,2n-1}(z,\overline{z})})\dznb\right)\\
&=& i\left(\dzj + a_j(z)\dzn+b_j(z)\dznb\right)\\
&& + \left( \tilde{L}_{2n-1,2j-1}(z,\overline{z}) +a_j(z)\tilde{L}_{2n-1,2n-1}(z,\overline{z}) +b_j(z)\tilde{L}_{2n-1,2n}(z,\overline{z})\right)\dzn\\
&& + \left(  \overline{\tilde{L}_{2n-1,2j}(z,\overline{z})}+a_j(z)\overline{\tilde{L}_{2n-1,2n}(z,\overline{z})}+b_j(z)\left( -2i+\overline{\tilde{L}_{2n-1,2n-1}(z,\overline{z})})\right)  \right) \dznb\\
&=& i X_j \text{ en remplaçant } a_j(z) \text{ et } b_j(z) \text{ par leurs valeurs données dans (\ref{a_i}) et (\ref{b_i})}\\
&&\text{ et en utilisant (\ref{J1}) et (\ref{J2}).} \\
\end{eqnarray*}
De plus,
\begin{eqnarray*}
J X_n &=& J \left( \dzn+b_n(z)\dznb \right)\\
&=&\left(i+\tilde{L}_{2n-1,2n-1}(z,\overline{z})\right)\dzn+\overline{\tilde{L}_{2n-1,2n}(z,\overline{z})}\dznb\\
&&+b_n(z)\left( \tilde{L}_{2n-1,2n}(z,\overline{z}) \dzn +\left(-i+\overline{\tilde{L}_{2n-1,2n-1}(z,\overline{z})})\dznb \right)\right) \\
&=&i\left(\dzn +b_n(z) \dznb\right)\\
&&+\left( \tilde{L}_{2n-1,2n-1}(z,\overline{z})+b_n(z) \tilde{L}_{2n-1,2n}(z,\overline{z})\right) \dzn\\
&&+\left(\overline{\tilde{L}_{2n-1,2n}(z,\overline{z})} +b_n(z)\left( -2i+\overline{\tilde{L}_{2n-1,2n-1}(z,\overline{z})}) \right)\right) \dznb\\
&=&iX_n \text{ en remplaçant }  b_n(z) \text{ par sa valeur donnée dans (\ref{b_n}) et en utilisant (\ref{J3}), (\ref{J4})}.\\
\end{eqnarray*}
$\square$

\subsection{Etude de l'espace tangent $J$-holomorphe $H^{1,0}_J\Gamma$.}
Nous nous intéressons maintenant à l'espace tangent $J$-holomorphe pour une hypersurface $\Gamma$ définie par $\Gamma=\{z\in\C^n,\,\rho(z)=0\}$, où $\rho$ est une fonction analytique réelle telle que $\rho(z)=\RE(z_n)+|z'|^2+o(|z|^2)$.\\
Soit $(L_1(z),\ldots,L_{n-1}(z))$ une base de $H^{1,0}_{J}\Gamma$ telle que, pour $i=1,\ldots,n-1$, $L_i(0)=\dzi$. Les champs $L_i$ sont des combinaisons linéaires des champs $X_j$, $j=1,\ldots,n$. On a donc :
\begin{eqnarray}
 L_i(z)=\dzi+\sum_{j=1}^{n}\alpha_{i,j}(z,\overline{z})\dfrac{\partial}{\partial z_j} + \beta_{i,n}(z,\overline{z})\dznb.
\end{eqnarray}
De plus, on a les développements limités suivants :
\begin{eqnarray*}
 \alpha_{i,j}(0)&=0, \, \text{pour} \,  j=1,\ldots,n-1\\
 \alpha_{i,n}(z,\overline{z})&=-(\tilde{b_i}+2\overline{z_i})+o(1),\\
\beta_{i,n}(z,\overline{z})&=\tilde{b_i}+o(1),\\
\end{eqnarray*}
où $\tilde{b_i}$ désigne le développement limité à l'ordre $1$ de $b_i$ (définition (\ref{b_i})).\\

Nous définissons un champ de vecteurs $T$ comme étant la projection du champ de vecteurs $[L_1,\overline{L_1}]$ dans $\C T\Gamma / H^{1,0}\Gamma\bigoplus\overline{H^{1,0}\Gamma}$. Nous avons alors $\C T\Gamma=H^{1,0}\Gamma\bigoplus\overline{H^{1,0}\Gamma}\bigoplus<T>$. Ainsi, $\{T,\, L_i,\,\overline{L}_i,\, i=1,\ldots , n-1\}$ est une base de $\C T \Gamma$, le complexifié de l'espace tangent de $\Gamma$.\\
Nous avons $T=i(\dzn-\dznb)$ordre 1 
Les règles de calcul ((\ref{rho1}), ... ) utilisées dans le cas modèle ne sont plus vérifiées. Nous pouvons cependant écrire :
\begin{eqnarray}
[L_j,\overline{L}_k]&=&\gamma_{j,\overline{k}}(z)T+\sum_{l=1}^{n}a_{j,k}^l(z)L_l+b_{j,k}^l(z)\overline {L_l},\label{rho2}\\
TL_k&=&L_kT+\gamma_k(z)T+\sum_{l=1}^{n}a_{k}^l(z)L_l+b_{k}^l(z) \overline {L_l}.\label{rho3}
\end{eqnarray} 
avec les développements limités suivants,
\begin{eqnarray*}
 \gamma_{j,\overline{k}}(z)&=&-2i\delta_{j,k}+\dfrac{1}{2}\left(\beta_j^k+\overline{\beta}_k^j\right)+ o(1),\\
a_{j,k}^l(0)&=&0,\\
b_{j,k}^l(0)&=&0,\\
\gamma_k(0)&=&0,\\
a_{k}^l(0)&=&0,\\
b_{k}^l(0)&=&0,\\
\end{eqnarray*}
On vérifie en particulier que l'on a encore, pour $j=1,\ldots,n-1$, $\gamma_{j,\overline{k}}(0)\neq 0$.

Soit $(Z_1(w),\ldots,Z_{n-1}(w))$ une base de $H^{1,0}_{J'}\Gamma'$, avec $Z_i=\dfrac{\partial}{\partial w_i}+\sum_{j=1}^{n}\alpha_{i,j}(w)\dfrac{\partial}{\partial w_j} + \beta_{i,n}(w)\dfrac{\partial}{\partial \overline{w_n}},\; i=1,\ldots , n-1$. Les équations de Cauchy-Riemann tangentielles vérifiées par la fonction $g$ sont les suivantes, pour $z\in\Gamma$ : 
\begin{eqnarray}
L_p \overline{g_j}&=0 \text{ pour }p=1,\ldots ,n-1 \text{ et } j=1,\ldots,n-1,\\
L_p \overline{g_n}&=\sum_{l=1}^{n-1}\dfrac{\beta_{l,n}(g)}{1+\varphi_l(g)}L_p g_l \text{ pour }p=1,\ldots ,n-1,\label{3.1bis'}\\
 L_p g_n&=\sum_{l=1}^{n-1}\dfrac{\alpha_{l,n}(g)}{1+\varphi_l(g)}L_p g_l, \text{ pour }p=1,\ldots ,n-1, \label{3.1ter'}\\
\dfrac{g_n+\overline{g_n}}{2}&+\sum_{j=1}^{n-1}g_j\overline{g_j}+\sum_{A,B}c'_{A,B}g^A\overline{g}^B=0.\label{3.2'}
\end{eqnarray}
où (\ref{3.2'}) est l'écriture de $ g(\Gamma)\subset \Gamma'$, avec $\rho'(w)=\RE(w_n)+|w'|²+\sum_{A,B}c'_{A,B}w^A\overline{w}^B$ et où $\varphi_l(w)=\sum_{m=1}^{n-1}\alpha_{m,l}(w)$.\\
Ces équations sont similaires aux équations (\ref{3.1}), (\ref{3.1bis}), (\ref{3.1ter}) et (\ref{3.2}) obtenues dans le cas modèle. Nous pouvons donc appliquer le procédé utilisé dans le cas des structures modèles. Les calculs sont plus techniques puisque les règles de calcul (\ref{rho2}) et (\ref{rho3}) entre les champs font intervenir des termes supplémentaires. Cependant, les termes supplémentaires qui apparaissent sont toujours de degré inférieur et ne constituent que des difficultés d'écriture. $\square$\

\subsection{Interprétation géométrique des structures presque complexes vérifiant la condition $(*)$}
 Soit $J$ une structure presque complexe sur $\C^n$ vérifiant la condition $(*)$, et soient $(X_1,\ldots,X_n)$ les champs $(1,0)$ tels que $X_j(0)=\dfrac{\partial}{\partial z_j}$. D'après (\ref{Xi}) et (\ref{Xn}), nous avons :
\begin{eqnarray*}
X_i&=&\dzi+\alpha_i(z)\dzn+\beta_i(z)\dznb,\; i=1,\ldots , n-1,\\
X_n&=&\dzn+\beta_n(z)\dznb,
\end{eqnarray*}
Le calcul des crochets de Lie donne, pour $j,k=1,\ldots , n-1$ :
\begin{eqnarray*}
[X_j,X_k]&=a_{j,k}(z)\dzn+b_{j,k}(z)\dznb,\\
\left[X_j,X_n\right]&=a_{j,n}(z)\dzn+b_{j,n}(z)\dznb,
\end{eqnarray*}
où, pour $j,k=1,\ldots,n-1$ :
\begin{eqnarray*}
a_{j,k}(z)&=&\dfrac{\partial \alpha_k}{\partial z_j}-\dfrac{\partial \alpha_j}{\partial z_k}+\alpha_j\dfrac{\partial \alpha_k}{\partial z_n}-\alpha_k\dfrac{\partial \alpha_j}{\partial z_n}+\beta_j\dfrac{\partial \alpha_k}{\partial \overline{z_n}}-\beta_k\dfrac{\partial \alpha_j}{\partial \overline{z_n}},\\
b_{j,k}(z)&=&\dfrac{\partial \beta_k}{\partial z_j}-\dfrac{\partial \beta_j}{\partial z_k}+\alpha_j\dfrac{\partial \beta_k}{\partial z_n}-\alpha_k\dfrac{\partial \beta_j}{\partial z_n}+\beta_j\dfrac{\partial \beta_k}{\partial \overline{z_n}}-\beta_k\dfrac{\partial \beta_j}{\partial \overline{z_n}},\\
a_{j,n}(z)&=&-\dfrac{\partial \alpha_j}{\partial z_n}+\alpha_j\dfrac{\partial \beta_n}{\partial z_n}-\beta_n\dfrac{\partial \alpha_j}{\partial \overline{z_n}},\\
b_{j,n}(z)&=&\dfrac{\partial \beta_n}{\partial z_j}-\dfrac{\partial \beta_j}{\partial z_n}+\beta_j\dfrac{\partial \beta_n}{\partial \overline{z_n}}-\beta_n\dfrac{\partial \beta_j}{\partial \overline{z_n}}.
\end{eqnarray*}
Nous avons donc :
\begin{lem}
Pour $j=1,\ldots,n-1$ et $k=1,\ldots,n$, les crochets de Lie $[X_j,X_k]$ vérifient $[X_j,X_k]=A_{j,k}X_n+B_{j,k}\overline{X}_n$, avec
\begin{eqnarray*}
A_{j,k}&=&\dfrac{a_{j,k}-\overline{\beta_n}b_{j,k}}{1-|\beta_n|^2},\\
B_{j,k}&=&\dfrac{b_{j,k}-\overline{\beta_n}a_{j,k}}{1-|\beta_n|^2}.
\end{eqnarray*}
 \end{lem}
Si $\Pi$ désigne la projection de $T\C^n$ sur $T\C^n / T^{1,0} \C^n$, on a $\Pi([X_j,X_k])=b_{j,k}\overline{X}_n$. Le terme $b_{j,k}\overline{X}_n$ précise donc le défaut d'intégrabilité de la structure : ce défaut est porté uniquement par la direction du champ de vecteurs $\overline{X_n}$.\\
Pour démontrer un résultat général concernant l'analyticité des applications CR, il est naturel de chercher à ''redresser`` les structures presque complexes, c'est-à-dire, de trouver un difféomorphisme local $\phi$ tel que $\phi_\star(J)$ soit une structure presque complexe vérifiant la condition $(*)$. L'analyse précédente indique que ceci n'est pas toujours possible. En effet, le fait que le défaut d'intégrabilité soit porté par une seule direction est stable par difféomorphisme. Une structure presque complexe qui ne vérifie pas cette condition ne pourra donc pas être redressée par un changement de coordonnées en une structure vérifiant la condition $(*)$. \\
Par exemple, dans le cas où $n=2$, soit $\C^2=(z_1,z_2)$, avec $z_j=x_j+iy_j$, $j=1,2$. Soit $J$ la structure presque complexe définie sur $\C^2$ par :
\begin{eqnarray*}J_\C=\begin{pmatrix}
i+a&b&0&0\\
\overline{b}&-i+\overline{a}&0&0\\
0&0&i+c&d\\
0&0&\overline{d}&-i+\overline{c}\\                                                                                                                                                                                                                                                                                                                                                                                                                                                                                                                                                                                                                                                                                                                                                                                                                                                                                                                                                                                                                                                                                                                                                                                                                                                                                                                                                                                                                                      \end{pmatrix},\end{eqnarray*} avec,\\
$\begin{array}{llll}
 a&=iy_2^2,&b&=y_2\sqrt{2+y_2^2},\\
c&=y_1\sqrt{2+y_1^2},&d&=iy_1^2.
\end{array}$\\
Les champs $(1,0)$ sont donnés par $X_1=\dfrac{\partial}{\partial z_1}+\alpha \dfrac{\partial}{\partial \overline{z_1}}$ et $X_2=\dfrac{\partial}{\partial z_2}+\beta\dfrac{\partial}{\partial \overline{z_2}}$, avec 
\begin{eqnarray*}
 \alpha=\dfrac{b}{2i-\overline{a}},&\beta=\dfrac{d}{2i-\overline{c}}.
\end{eqnarray*}
On a 
\begin{eqnarray*}
 [X_1,X_2]&=&\dfrac{\partial \beta}{\partial z_1}\dfrac{\partial}{\partial \overline{z_2}}-\dfrac{\partial \alpha}{\partial z_2}\dfrac{\partial}{\partial \overline{z_1}},\\
&=&\dfrac{\overline{\alpha}}{1-|\alpha|^2}\dfrac{\partial \alpha}{\partial z_2}X_1-\dfrac{1}{1-|\alpha|^2}\dfrac{\partial \alpha}{\partial 
z_2}\overline{X_1}-\dfrac{\overline{\beta}}{1-|\beta|^2}\dfrac{\partial \beta}{\partial 
z_1}X_2+\dfrac{1}{1-|\beta|^2}\dfrac{\partial \beta}{\partial z_1}\overline{X_2}.
\end{eqnarray*}
D'où,
\begin{eqnarray*}
\Pi([X_1,X_2])&=-\dfrac{1}{1-|\alpha|^2}\dfrac{\partial \alpha}{\partial 
z_2}\overline{X_1}+\dfrac{1}{1-|\beta|^2}\dfrac{\partial \beta}{\partial z_1}\overline{X_2}.
\end{eqnarray*}
Or,
\begin{eqnarray*}
-\dfrac{1}{1-|\alpha|^2}\dfrac{\partial \alpha}{\partial 
z_2}&=\dfrac{1}{2\sqrt{2+y_2^2}},\\
\dfrac{1}{1-|\beta|^2}\dfrac{\partial \beta}{\partial z_1}&=\dfrac{-1}{2\sqrt{2+y_1^2}}.\\
\end{eqnarray*}
La projection $\Pi([X_1,X_2])$ n'est pas portée par une direction constante. On ne peut donc pas redresser la structure $J$ en une structure vérifiant la condition $(*)$ par un difféomorphisme local.

\providecommand{\bysame}{\leavevmode\hbox to3em{\hrulefill}\thinspace}
\providecommand{\MR}{\relax\ifhmode\unskip\space\fi MR }
\providecommand{\MRhref}[2]{%
  \href{http://www.ams.org/mathscinet-getitem?mr=#1}{#2}
}
\providecommand{\href}[2]{#2}

\vskip 0,5cm
{\small
\noindent Marianne Peyron\\
(1) UJF-Grenoble 1, Institut Fourier, Grenoble, F-38402, France\\
(2) CNRS UMR5582, Institut Fourier, Grenoble, F-38041, France\\
{\sl E-mail address} : marianne.peyron@ujf-grenoble.fr
}
\end{document}